\documentclass[ijoc,nonblindrev]{informs3} % current default for manuscript submission

\OneAndAHalfSpacedXII % current default line spacing
\usepackage{natbib}
 \bibpunct[, ]{(}{)}{,}{a}{}{,}%
 \def\bibfont{\small}%
\TheoremsNumberedThrough     % Preferred (Theorem 1, Lemma 1, Theorem 2)
\EquationsNumberedThrough    % Default: (1), (2), ...
\newcommand{\textwrap}{\parbox[t]{4.5in}}
\newcommand{\textwrapbig}{\parbox[t]{5.5in}}

\def \Ehat{\hat E}

\def \Rhat{\hat R}

\def \hbar{\bar h}

\def \Vbar{\overline V}

\def \Bcal{{\mathcal B}}

\def \Fcal{{\mathcal F}}
\def \Gcal{{\mathcal G}}
\def \Hcal{{\mathcal H}}
\def \Ical{{\mathcal I}}

\def \Qcal{{\mathcal Q}}

\def \Xcal{{\mathcal X}}
\def \Ycal{{\mathcal Y}}

\usepackage{times}                 % uncomment these lines
\usepackage{amsmath}
\usepackage{nccmath}
\usepackage{amssymb,color, graphicx}
\usepackage{multirow}
\usepackage{epsfig}
\usepackage{bbm}
\usepackage{algorithmic}
\usepackage{algorithm}
\usepackage{pgfplots}
\usepackage{tikz}
\usepackage[inline]{enumitem}
\usepackage{hyperref}
\usepackage{slashbox}
\usepackage{comment}
\usetikzlibrary{pgfplots.groupplots}
\usetikzlibrary{external}
\tikzset{external/mode=graphics if exists}
\tikzset{
    png export/.style={
        external/system call/.add={}%
        {; convert -density 300 -transparent white "\image.pdf" "\image.png"}
    }
}
\tikzset{external/system call={lualatex
     \tikzexternalcheckshellescape -halt-on-error
     -interaction=batchmode -jobname "\image" "\texsource"}}

\usepgfplotslibrary{external}
\usepgfplotslibrary{dateplot}

\newcommand{\E}{\mathbbm{E}}

\newcommand{\bn}{\begin{eqnarray}}
\newcommand{\en}{\end{eqnarray}}
\newcommand{\bns}{\begin{eqnarray*}}
\newcommand{\ens}{\end{eqnarray*}}

\newcommand{\F}{\mathcal{F}}

\newcommand{\Ic}{\mathcal{I}}

\newcommand{\Nc}{\mathcal{N}}
\newcommand{\Ec}{\mathcal{E}}

\newcommand{\Xc}{\mathcal{X}}
\newcommand{\Gc}{\mathcal{G}}
\newcommand{\R}{\mathbbm{R}}

\renewcommand{\Delta}{\varDelta}
\renewcommand{\Pi}{\varPi}
\newcommand{\Pb}{\mathbb{P}}

\def \Vbar{\overline V}

\definecolor{pinegreen}{rgb}{0.1,0.5,0.2}

\begin{document}
%%%%%%%%%%%%%%%%

% Outcomment only when entries are known. Otherwise leave as is and 
%   default values will be used.
%\setcounter{page}{1}
%\VOLUME{00}%
%\NO{0}%
%\MONTH{Xxxxx}% (month or a similar seasonal id)
%\YEAR{0000}% e.g., 2005
%\FIRSTPAGE{000}%
%\LASTPAGE{000}%
%\SHORTYEAR{00}% shortened year (two-digit)
%\ISSUE{0000} %
%\LONGFIRSTPAGE{0001} %
%\DOI{10.1287/xxxx.0000.0000}%

% Author's names for the running heads
% Sample depending on the number of authors;
% \RUNAUTHOR{Jones}
% \RUNAUTHOR{Jones and Wilson}
\RUNAUTHOR{Asamov, Salas and Powell}
% \RUNAUTHOR{Jones et al.} % for four or more authors
% Enter authors following the given pattern:
%\RUNAUTHOR{}

% Title or shortened title suitable for running heads. Sample:
% \RUNTITLE{Bundling Information Goods of Decreasing Value}
% Enter the (shortened) title:
\RUNTITLE{SDDP vs. ADP}

% Full title. Sample:
% \TITLE{Bundling Information Goods of Decreasing Value}
% Enter the full title:
%\TITLE{}
\TITLE{SDDP vs. ADP: The Effect of Dimensionality in Multistage Stochastic Optimization for Grid Level Energy Storage}

% \date{Received: date / Accepted: date}
\date{}

\ARTICLEAUTHORS{%
\AUTHOR{Tsvetan Asamov, Daniel F. Salas and Warren B. Powell}
\AFF{Department of Operations Research and Financial Engineering, Princeton University\\
\EMAIL{\{tasamov, dsalas, powell\}@princeton.edu}}
% \AUTHOR{Warren B. Powell}
% \AFF{Department of Operations Research and Financial Engineering, Princeton University, \EMAIL{powell@princeton.edu}}
% Enter all authors
} % end of the block

\ABSTRACT{%
There has been widespread interest in the use of grid-level storage to handle the variability from increasing penetrations of wind and solar energy.  This problem setting requires optimizing energy storage and release decisions for anywhere from a half-dozen, to potentially hundreds of storage devices spread around the grid as new technologies evolve.  We approach this problem using two competing algorithmic strategies.  The first, developed within the stochastic programming
literature, is stochastic dual dynamic programming (SDDP) which uses Benders decomposition to create a multidimensional value function approximations, which have been widely used to manage hydro reservoirs.  The second approach, which has evolved using the language of approximate dynamic programming, uses separable, piecewise linear value function approximations, a method which has been successfully applied to high-dimensional fleet management problems.  This paper brings these two
approaches together using a common notational system, and contrasts the algorithmic strategies (which are both a form of approximate dynamic programming) used by each approach.  The methods are then subjected to rigorous testing using the context of optimizing grid level storage.
% , which provides a new setting characterized by potentially high-dimensional resource vectors, temporally fine-grained processes, and diverse information processes.
}%

% Sample 
%\KEYWORDS{deterministic inventory theory; infinite linear programming duality; 
%  existence of optimal policies; semi-Markov decision process; cyclic schedule}

% Fill in data. If unknown, outcomment the field
\HISTORY{}

\maketitle

\section{Introduction}
There is global interest in increasing the generation of electricity from renewables to meet the pressure to reduce our carbon footprint.  However, wind and solar energy cannot be directly controlled (they are not ``dispatchable'' in the parlance of the energy systems community), and in addition to their predictable variability (e.g. the rising and setting of the sun), they introduce a high level of uncertainty.  Over the years, grid operators have developed a sophisticated planning process to plan power needs in the presence of modest levels of uncertainty, due primarily to changes in weather.  There is a growing consensus that storage will be needed to smooth the variations introduced by wind and solar, which requires making decisions about when and where to charge and discharge batteries (or other storage devices).  Storage devices will need to be managed in a way that captures both the effects of grid congestion, as well as the real-time control of generators.

The process of managing energy generation and transmission by grid operators in the U.S. consists primarily of three steps: 1) the day-ahead unit-commitment problem, which determines which steam-generating units will be turned on and off (and when); 2) intermediate-term planning (typically every 15 to 30 minutes) which determines which gas turbines will be turned on and off, and 3) real-time economic dispatch (every 5 minutes), which is the process of increasing/decreasing (``ramping'')
the output of generators (both steam and gas turbines) in response to current conditions.  To handle unexpected variations (due to weather and load variations which are caused in part by the sun), grid operators plan reserve capacity, typically in two forms: spinning reserve, which can be tapped immediately, or nonspinning reserve, which is usually gas turbines that can be brought online in a few minutes.

This process handles unexpected variations, but not at the level that would be experienced at the high penetrations of renewables that are being targeted by state renewable portfolio standards in the U.S., as well as national targets being set by nations around the world.  As grid operators draw close to 20 percent renewables from wind and solar (in Germany it is over 30 percent), it is widely anticipated that grid-level storage will be needed to help smooth the variations to meet grid capacity constraints and to handle the gaps between the rate of change from wind and solar, and the ramping ability of online generators.

Energy storage is an emerging technology at this time.  A major grid such as that operated by PJM may have only 5-10 storage devices (batteries and pumped hydro) but the number is growing quickly as the technology improves and costs come down.  In addition to the possibility of batteries being installed in individual communities to help with outages, vehicle-to-grid technology already exists that allows an energy aggregator to treat a few hundred cars as a single storage ``device.''  It is easy to envision a grid with hundreds of large storage devices, and perhaps thousands of smaller ones.

This paper addresses the algorithmic challenge of simultaneously optimizing decisions to charge and discharge energy in a network of grid-connected batteries, while also optimizing ramping decisions at each generator that is currently operating.  These decisions have to also respect transmission constraints.  We wish to test this logic at different levels of investment in solar generation capacity.  To do this, we first use a model called SMART-ISO which simulates the day-ahead, intermediate and real-time planning processes at PJM; this model was carefully calibrated against historical performance at PJM and was shown to accurately replicate network behavior \cite{maowit2015b}.  However, SMART-ISO, which takes approximately 2-3 hours to simulate a single week, is unable to optimize storage.  For this reason, we use the unit commitment decisions from SMART-ISO (which adapts to any level of solar investment we would like to simulate), but then use a model that optimizes only the ramping decisions along with grid congestion, while also optimizing charge and discharge decisions for storage devices.

A method for optimizing grid-level storage was recently proposed in \cite{Salas2015} based on approximate dynamic programming.  This method uses separable, piecewise linear value functions to capture the value stored in each device around the grid, at 5-minute increments.  The method was shown to produce near-optimal solutions for deterministic problems, which was the only benchmark available at the time.  This leaves open the very real question of how well this methodology would work under more realistic stochastic conditions.

Far more popular in the stochastic programming community has been the use of Benders decomposition, primarily in the context of a methodology known as stochastic dual dynamic programming (SDDP), which was first proposed in \cite{Pinto1991} for the management of water reservoirs.  This has produced a flurry of follow-on papers (\cite{Mo2001}, \cite{Shapiro2011}, \cite{LoWoMi2013}) focusing primarily on applications in the planning of hydroelectric power, as well as
considerable theoretical interest (\cite{Linowsky2005}, \cite{Philpott2008}, \cite{ShDeRu2014}, \cite{Girardeau2012a}).  SDDP has generally made the assumption of ``intertemporal independence'' which is that new information becoming available at time $t$ does not depend on the history of the process.  Benders has been used in conjunction with scenario trees that capture the history of the information process (\cite{Birge1985}, \cite{Sen2014a}), but this work has not proven to
be computationally tractable.  Further, even with the assumption of intertemporal independence, SDDP can only be applied to a sampled approximation.  Researchers have argued that the use of a sampled model produces small errors (see e.g. \cite{ShDeRu2014}), but we are unaware of any research evaluating errors in specific decisions, such as the congestion on a transmission line that might guide hundreds of millions of dollars in new investment.  Finally, it is generally well known that Benders
struggles with high dimensional problems, where ``high'' might mean more than 20 storage devices.  However, we are unaware of any formal study of errors that arise when using Benders (or separable approximations) as a function of the dimensionality of the resource vector. This is a question we address in this paper.  

Recent research has mitigated some of the limitations of classical SDDP.  \cite{Asamov2015} presented a version of Benders decomposition with a new regularization strategy designed for multistage problems (regularization has long been recognized as a powerful technique for Benders, but this is the first extension to multiperiod problems).  Further, this work also considers a version that exhibits a first-order Markov information process.

This paper, then, uses the setting of grid level storage on the PJM grid to compare two algorithmic strategies:
\begin{itemize}
\item[1)] Regularized Benders decomposition for multiperiod (and multistage) problems (SDDP), with independent and first-order Markov information processes.  This method can only be run on a sampled information process.
\item[2)] Approximate dynamic programming using separable, piecewise linear value function approximations (ADP-SPWL).  This method can also be run using independent and first-order Markov information processes, but is not limited to a sampled version of the problem.
\end{itemize}
These experiments will provide the first serious benchmark for both algorithmic strategies.  Benders decomposition has long enjoyed bounds, but for our problem, we show that these bounds are not very tight, and further provide little insight into the accuracy of individual decisions which may affect investment decisions in the grid.  ADP using piecewise linear value functions has been studied over the years in transportation applications (e.g. \cite{Topaloglu2006}), but in the past the only benchmarking has been against deterministic solutions.  We would also argue that up to now, Benders decomposition has not faced serious algorithmic competition.

This paper makes the following contributions: 1) We present two algorithmic strategies, SDDP (from stochastic programming) and ADP-SPWL (from approximate dynamic programming) and show that these are both forms of approximate dynamic programming where the primary difference is how the value functions are being approximated (multidimensional Benders cuts vs. separable, piecewise linear approximations).  We then highlight other differences that result specifically because of
the nature of the two approximation strategies.  2) We then use a model of the PJM power grid and real-time energy generation to optimize across storage portfolios ranging from 5 to 100 batteries, which allows us to test the quality of the solution from each algorithmic strategy as a function of the dimensionality of the resource state variable.  This appears to be the first comprehensive evaluation of SDDP over a wide range of dimensions, and the first formal evaluation of the ADP
approach on a multidimensional stochastic problem, using the context of energy storage that introduces much higher dimensionality than has appeared in other experiments.  In addition, our problem setting involves more time periods than are typically considered (288, representing 24 hours in 5-minute increments), for a problem that is highly time-dependent.  3) We use the ADP strategy to provide the first evaluation of errors introduced in SDDP by solving a sampled model, focusing on the quality of individual decisions.
% 4) We close with a study of high levels of investment in solar capacity to study how the grid behaves when simultaneously optimizing across storage and real-time ramping of generations.

The paper is organized as follows. Section \ref{sec:gridpapermodel} provides a mathematical model of the grid-level storage problem. Section \ref{sec:gridpaperadp} presents canonical models of multistage stochastic programs, and dynamic programs and demonstrates that SDDP and ADP share the same structure, with just minor (but important) differences in approximation strategies.  Section \ref{sec:gridpaperexperiments} provides a thorough set of comparisons between SDDP and ADP-SPWL in an
energy storage setting, where we can vary the dimension of the resource state variable from 5 to 100, representing the first serious test of both algorithmic strategies over a wide range of dimensionalities.  % Then, section \ref{sec:gridlevelsolar} reports on a series of experiments to study the behavior of high-levels of solar penetration on the PJM grid.
Section \ref{sec:gridpaperconclusions} concludes the paper.

\section{Mathematical model of grid-level storage}
\label{sec:gridpapermodel}
Our intent is to study storage in the presence of high levels of energy from off-shore wind farms, derived as a part of a larger study of off-shore wind.  The careful modeling of offshore wind is given in \cite{Archer2015}.  We then used a large-scale model of the PJM grid and energy markets called SMART-ISO to plan the correct level of energy from slow and fast fossil generating plants for a given level of wind penetration (see \cite{maowit2015b} for a thorough description of this study).  SMART-ISO models the nested planning process consisting of day-ahead, intermediate (roughly hour-ahead) and real-time planning, closely matching PJM's actual planning process.  SMART-ISO carefully replicates the uncertainty in each planning step.  For example, forecast errors in planning day-ahead and hour-ahead energy from wind are based on actual errors from the forecasts that PJM uses for its planning of its own (onshore) wind farms.

SMART-ISO is a large-scale simulator of the unit commitment process that is unable to optimize storage.  For this reason, we would run SMART-ISO assuming a particular investment in wind generation capacity.
In our setup, SMART-ISO models power plants of various types with a total of 825 generators and a combined maximum capacity of 129,638 MW. We consider 396 gas turbines (23,309 MW), 50 combined cycle generators (21,248 MW), 264 steam generators (73,374 MW), 31 nuclear reactors (31, 086 MW), and 84 conventional hydro power generators (2,217 MW).
This model would then determine which generators were on or off at any given point in time. We then passed these on/off decisions to the storage model which would then optimize energy storage decisions in the presence of grid congestion.  Optionally, our storage model is also allowed to optimize ramping decisions of fossil generators (without turning them on or off), although we had to turn off this feature for some of our algorithmic testing.  

Below, we present our mathematical model of the grid level storage problem, spanning storage, transmission, energy generation from fossils, and the exogenous generation of energy from wind.

\subsection{Problem Description}
\label{s:storage}
We consider the grid of PJM Interconnection (or simply PJM), a large independent regional transmission operator serving the mid-Atlantic states out to Chicago. The PJM network comprises more than 9,000 buses and 14,000 transmission lines. A map of the the geographical territory that is served by PJM is shown in Figure 1 in the online supplement. %\ref{f:map}.
%Daily operation schedules are determined by solving a unit commitment problem consisting of three parts: the day-ahead, the hour-ahead and the real-time economic dispatch models. The day-ahead and hour-ahead solutions determine the on/off status of power generators for respectively the next day and next hour. Additionally, the generation levels of active units are adjusted by the real-time economic dispatch decisions in order to satisfy the observed demand. The control of multiple storage devices in a distributed energy system is a challenging task that depends on a variety of factors such as the location of each device, the presence of transmission line congestion, and a variety of other constraints. Furthermore, the complexity of the problem increases significantly as the number of storage devices grows. Our inputs comprise real data for the parameters of the power generators such as the operational characteristics, variable cost, and location in the grid. Our goal is to solve a grid-scale model for the economic dispatch problem in the presence of renewables.

\textbf{The Grid}\\
We model the PJM grid as a graph $\Gc = (\Nc, \Ec)$. The set of nodes $\Nc$ correspond to buses in the grid and the set of edges $\Ec$ correspond to transmission lines. More specifically, each edge $(n_i, n_j)\in\Ec$ corresponds to the transmission line connecting bus $n_i\in\Nc$ to bus $n_j\in\Nc$. Furthermore, we use the DC power flow approximation to model the power flow between buses in the grid.
% \begin{figure}
% \begin{center}
% \includegraphics[width=4in]{Figs/PJMgrid.jpg}
% \end{center}
% \caption{The PJM grid.}
% \label{f:map}
% \end{figure}

\textbf{Parameters}\\
% We consider a total of 825 generators of various types such as gas turbines, combined cycle, steam, nuclear, and hydro power generators.
% The combined maximum capacity of all power generators in the model is 129,638 MW.
The set of parameters for each generator includes its power capacity and generation cost. % , as well as maximum ramp-up and ramp-down rates.
Each storage device is characterized by its minimum and maximum energy capacity, its charging and discharging efficiency, and its variable storage cost. To present a formal description of our model, we introduce the following notation:
\bns
\Gcal &=& \textwrap{The set of fossil generators,} \nonumber \\
\Bcal &=& \textwrap{The set of storage devices,}\nonumber \\
\Qcal &=& \textwrap{The set of wind farms,}\nonumber \\
\kappa^l_g, \kappa^u_g &=& \textwrap{The minimum and maximum power capacities for electricity generator $g\in\Gcal$,}\nonumber \\
\kappa^l_b, \kappa^u_b &=& \textwrap{The minimum and maximum energy capacities for storage devices $b\in\Bcal$,}\nonumber \\
\eta^+_b, \eta^-_b &=& \textwrap{The charging and discharging multipliers (efficiencies) for each storage device $b\in\Bcal$,}\nonumber \\
c^G_{t,g} &=& \textwrap{The vector of variable generation cost in \$/MWh for the set of generators $g\in\Gcal$ at time $t$,}\nonumber \\
c^B_{t,b} &=& \textwrap{The vector of variable storage cost in \$/MWh for the set of storage devices $b\in\Bcal$ at time $t$,}\nonumber \\
d_t&=& \textwrap{The vector of electricity demands (loads), in MW, for each node at time $t=0,\dots, T$. We assume that electricity demand evolves deterministically,}\nonumber \\
%Z^G_{tg} &=& \textwrap{The binary variable indicating whether generator $g\in\Gcal$ is on or off (determined by SMART-ISO) during time period $t$,}\nonumber \\
\Ycal &=& \textwrap{The closed and convex set describing the model for the DC power flow in the grid following Kirchhoff's laws. It takes into account the structure of the electrical grid, capacities of the transmission lines and bounds on phase angles (similar to voltage limits in AC power flow)}. \nonumber
\ens

\textbf{State variables}\\
We let $S_t$ be the (pre-decision) state capturing all the information we need at time $t$ to model the system from time $t$ onward, and we let $S^x_t$ be the post-decision state, which is the information we need after we have made a decision $x_t$ at time $t$.  The evolution of states, decisions and information is described by time:
$$S_0\xrightarrow{x_0} S^x_0\xrightarrow{\omega_1} S_1\xrightarrow{x_1} S^x_1\xrightarrow{\omega_2}  \ldots \xrightarrow{\omega_T} S_T \xrightarrow{x_T} S^x_T.$$
%\begin{definition}
%The (pre--decision) state $S_t$ of the system at time $t\geq 1$ is all the information that is necessary and sufficient to make a decision at time $t$, and model the computation of costs, constraints, and transitions from time $t$ onward.
%\end{definition}

Formally, we define the elements of the state of the system $S_t = (R_t, I_t)$ as consisting of resource state variables $R_t$ and (exogenous) information state variables $I_t$.
%The resource state variables are given by
The resource state $R_t$ of the system at time $t$ is a real vector that can be partitioned into subvectors as $R_t = (R^B_{t}, R^G_{t})$ where:
\bns
% R_t &=& \textwrap{The resource state $R_t$ of the system at time $t$ is a real vector that can be partitioned into subvectors as $R_t = (R^B_{t}, R^G_{t})$ where:}\nonumber \\
R^B_{t,b} &\textwrap{is the amount of energy available in storage device $b\in\Bcal$ at time $t$,}\nonumber \\
R^G_{t,g} &\textwrap{is the level of power output from generators at time $t$. Since generators have up- and down- ramping limits, the power output at time $t$ affects the range of feasible power output levels at time $t+1$.}\nonumber 
\ens
% We let $R^B_t$ be the vector giving the storage levels in our batteries at time $t$, and $R^G_t$ be the power output level of operating generators determined by the vector $Z^G$ defined below.  The resource state is then given by
% \bns
% R_t = (R^B_t,R^G_t),
% \ens
Please note that the power output levels $R^G_t$ obey the on/off generator status determined by SMART-ISO, which is given by the vector $Z^G$ defined below. 

% \begin{remark}
% The dimension of the post-decision resources state is $|R^x_t| = |R_{t,M}^x| + |R_{t,F}^x|$.
% \end{remark}
The information state $I_t$ includes the following variables: 
\bns
% I_t &=& \textwrap{The information state $I_t$ which includes:}\nonumber \\
E^W_t &=& \textwrap{The energy from wind at time $t$,}\nonumber \\
L_{ti} &=& \textwrap{The load at time $t$ at node $i$,}\nonumber \\
Z^G_{tg} &=& \textwrap{The binary variable indicating whether generator $g\in\Gcal$ is on or off (determined by SMART-ISO) during time period $t$.}\nonumber
\ens
We let $L_t$ be the vector of loads and $Z^G_t$ be the vector indicating which generators are turned on or off.  Thus our information state is a the following vector 
\bns
I_t = (E^W_t,L_t,Z^G_t), \ens
where $Z^G_t$ is provided exogenously, and $R^G_t = 0$, if $Z^G_t = 0$ (the generator output is set to 0 when it is turned off).

We also represent the post-decision state of the system as $S_t^x = (R_t^x, I_t^x)$ which is the state of our system immediately after a decision is made.

The post-decision resource state $R^x_t = (R^{B,x}_t, R^{G,x}_t)$ is given by the post-decision amount of energy available in storage devices,
and the post-decision level of power output from generators $R^{G,x}_t$. Using the decision notation introduced below, we define the post-decision resource state as the pre-decision resource state $R^B_t$ adjusted for battery inflows and outflows, 
\bns
R^{B,x}_{t, b}  = R^B_{t, b} + \eta^-_{b} x^{B-}_{t,b} - \eta_t^+  x^{B+}_{t,b}.
\ens
When modeling batteries, there is no exogenous change to the resource level, which means that the pre-decision state is given by
\bns
R^B_{t+1} = R^{B,x}_t.
\ens
If we model a water reservoir, we could account for exogenous rainfall (say, to reservoir $b$) using
\bns
R^B_{t+1,b} = R^{B,x}_{t,b} + \Rhat_{t+1,b}
\ens
where $\Rhat_{t+1,b}$ would be the stochastic rainfall occurring between $t$ and $t+1$.

%\begin{remark}
%The given definitions are sufficient for most practical purposes but might be considered a bit too informal for theoretical work. For a more formal approach, please see the extended version presented in \cite{asamov2015regularized}.
%\end{remark}

\textbf{Decisions}\\
% At each time period $t$, we need to determine the amount of power injected or withdrawn from the grid at each storage device and each power generator.
% We manage the dynamic allocation of energy over a time horizon of length $T$.
Decisions are made at each time step $t\in\{0,1,\dots, T\}$ to determine the energy flow between the electric grid, wind farms and storage devices subject to meeting electricity demand. Energy from wind farms can be used to satisfy the current demand, or it can be transmitted directly into the storage devices. At any time period $t$, energy available in storage devices can be sold to satisfy the grid demand. Furthermore, energy can also be bought from the grid and transferred into storage for later use.

We partition the decision vector $x_t\in\Xc_t$ as 
\bns
x_t = \begin{pmatrix}x^{G+}_{t}\\ x^{B+}_{t}\\ x^{B-}_{t}\\ y_t\end{pmatrix},
\ens
where $x^{G+}_{t}\in\R^{|\Gcal|}$, while $x^{B+}_{t}, x^{B-}_{t}\in\R^{|\Bcal|}$ and $y_t\in Y$. The vectors $x^{B+}_{t}, x^{G+}_{t,g}$ denote the respective amounts of power injected into the grid by storage devices and generators, and $x^{B-}_{t}$ denotes the power withdrawn from the grid by storage devices at time $t$.\\
Decisions have to reflect a number of constraints, which we describe next.
\begin{itemize}
\item When a generator $g\in \Gcal$ first comes online $\big((Z^G_{t,g} =1) \cap (t = 0 \cup Z^G_{t-1,g} = 0)\big)$, its power output is set to its minimum power capacity. Thus we have the initial generation constraints:
\begin{equation}
\label{eq:15}
\begin{split}
x^{G+}_{0,g} &= \kappa_g^l,~ \mbox{ for } g\in \Gcal ~\mbox{ such that } Z^G_{0,g}=1.\\
x^{G+}_{t,g} &= \kappa_g^l,~ \mbox{ for } t = 1,\dots, T,~\mbox{ and } g\in \Gcal ~\mbox{ such that } (Z^G_{t-1,g} = 0~\cap~ Z^G_{t,g}=1).
\end{split}
\end{equation}
%Power generator capacity constraints.
% \item Generator ramp-up and ramp-down constraints: Since generators have ramping limits the power output at time period $t$ affects the feasible range for the next time period.
% \item The level of output from power generators ().
\item The output of any active power generator is bounded by its minimum and maximum capacity. Furthermore, inactive generators (that is, where $Z^G_{t,g} = 0$) must have zero output.  Hence we have the following capacity constraints for the set of power generators:
\begin{equation}
\label{eq:16}
\kappa^l_g  Z^G_{t,g} \leq x^{G+}_{t,g}\leq \kappa^u_g  Z^G_{t,g},~\mbox{ for } t=0,\dots, T.
\end{equation}
% Storage capacity constraints for batteries:
% \begin{equation}
% \begin{split}
% \kappa^l_M \leq x_{t,M}^+\leq \kappa^u_M \\
% \kappa^l_M \leq x_{t,M}^-\leq \kappa^u_M
% \end{split}
% \end{equation}
\item
We let the vector $p_t\in \R^{|\Nc|}$ represent the nodal power generation at each node in $\Nc$. Given a node $n_i\in\Nc$,
we denote with $\Gcal(n_i)$, $\Bcal(n_i)$, and $\Hcal(n_i)$ the sets of respectively generators, storage devices, and wind farms that map to node $n_i$. Hence, the $i-$th component of $p_t$ is
\begin{equation}
\label{eq:17}
p_{t,i} = \sum_{g\in \Gcal(n_i)} x^+_{t,g} + \sum_{b\in \Bcal(n_i)} \big(x^+_{t,b} - x_{t,b}^- \big) + \sum_{q\in \Qcal(n_i)}E^W_{t,q},~ t=0,\dots,T,~ i=1,\dots, |N|.
\end{equation}
\item Since $y_{t,i}$ denotes the amount of power arriving at node $n_i\in N$ at time $t$, we can write the electricity demand constraints as
\begin{equation}
\label{eq:18}
y_t + p_t = d_t,~\mbox{ for } t=0,\dots, T.
\end{equation}
\item Furthermore, we impose flow conservation constraints for storage devices $b\in\Bcal, t=0,\ldots,T$:
\begin{equation}
\label{eq:19}
\kappa^l_b \leq R^B_{t,b} + \eta^- x^{B-}_{t,b} - {\eta^+} x^{B+}_{t,b} \leq \kappa^u_b.
\end{equation}
\end{itemize}

We let $\Xcal_t,~t=0,\ldots, T$ be the feasible region defined by the constraints (\ref{eq:15})-(\ref{eq:19}).

In the context of energy storage, complete recourse is an obvious property since we can always choose not to use the storage devices. This is similar to many other real-world applications where a default decision is always available.
% This is one of these concepts invented by mathematicians that never seems to arise in practice (). I have never seen a real-world application that needs feasibility cuts.

In our model, we do not impose ramping constraints for both power generators and the distributed storage devices.  While incorporating ramping constraints is relatively easy within the algorithm, it means that we need to include the current output level of each generator in our state variable, which complicates our ability to prove tight bounds.

\textbf{The exogenous information process}\\
The only exogenous information processes we consider in this paper is the stochastic change in energy from wind which we model using
\bns
\Ehat^W_t &=& \textwrap{The change in energy from wind between $t-1$ and $t$,} \\
\ens
Energy from wind can be modeled using a rolling forecast of wind provided by an exogenous source (this is how it is done at PJM).  These forecasts are provided as exogenous information in the form of a vector
\bns
f^W_{tt'} = \textwrapbig{The forecast of wind at time $t'$ using the information we have at time $t$.}
\ens
Using this notation, $E^W_t = f^W_{tt}$, and the exogenous change between $t$ and $t+1$ would be
\bns
\Ehat^W_{t+1} &=& f^W_{t+1,t+1} - f^W_{t,t+1}.
\ens
\begin{remark}
In some problems, additional (stochastic) renewable sources such as solar power might be present. In that case we can use similar notatation. For example, we can denote with
\bns
\Ehat^S_t &=& \textwrap{The change in energy from solar between $t-1$ and $t$.}
\ens
\end{remark}

In our numerical work, we view $L_t$ and $Z^G_t$ as exogenous information arriving to the system over time.  However, we do not model these as exogenous processes, which means we can treat the vectors $(L_t,Z^G_t)$, $t=0,\ldots T$ as {\it latent variables}.

Our wind modeling was derived from a study of off-shore wind reported in \cite{Archer2015}, which combined a base set of forecasts from a meteorological model called Weather, Research and Forecasting (WRF), and a stochastic model of errors in forecasts derived from historical data (and forecasts) provided by PJM.  Solar data was derived from actual solar energy (in 5-minute increments) from 23 solar farms operated by PSE\&G, a large utility based on New Jersey.  We factored up wind and solar data to test our algorithms at high levels of each source of energy.

The model is driven by two other deterministic (but time varying) processes.  These are
\bns
L_{ti}&=\textwrap{ the load (in MW) at time $t$ at node $i$,}\\
Z^G_{tg} &=
\begin{cases}
1,~ \textwrap{ if generator $g$ is on at time $t$},\\
0,~ \textwrap{ otherwise.}
\end{cases}
\ens
The loads $L_t$ are those served by PJM in 2010 (in 5-minute increments).  The indicator variables $Z^G_{tg}$ are determined by a unit commitment simulator called SMART-ISO which optimizes fossil generation based on forecasted and actual levels of energy from wind.  However, our grid model controls the output level of generators that are already turned on.  This way, we trade off generating energy from fossil generators, and storing energy in our grid-level storage devices.

We could introduce other sources of uncertainty such as variations in loads, but our plan is to model high penetrations of wind. At these levels, the uncertainty from wind forecasts is much higher than the uncertainty from other sources.  We do not consider outages due to the failure of generators or transmission lines.  This model allows us to focus on using storage purely to handle the variability due to wind or solar.

Even though $L_t$ and $Z^G_t$ are deterministic processes, we model them as if they are revealed over time.  For this reason, we let $W_t = (E^W_t,L_t,Z^G_t)$ be our exogenous information process.
More formally, we are given a probability space $(\Omega, \Fcal, \Pb)$ with a sigma algebra $\Fcal$, and a filtration $\{\emptyset, \Omega\} = \Fcal_1\subset\Fcal_2\subset\dots\subset\Fcal_T=\Fcal$. The stochastic process $\{W_t\}_{t=1}^T$ is adapted to $\{\Fcal_t\}_{t=1}^T$, and the sets of possible realizations of $W_t$ are denoted with $\Omega_t$. Those correspond to nested partitions of $\Omega$ that are given by $\{\Fcal_t\}_{t=1}^T$. 

% We let $\omega\in\Omega$ represent a sample realization of $W_1, W_2, \ldots, W_T$.  Let $\F$ be the sigma-algebra on $\Omega$, and let $\Pcal$ be a probability measure on $(\Omega,\Fcal)$, giving us a probability space $(\Omega,\Fcal,\Pcal)$.  We further let $\Fcal_t$ be the sigma-algebra generated by $W_1, \ldots, W_t$.  We adopt the convention that any variable indexed by $t$ is $\Fcal_t$-measurable.

%More formally, given a probability space $(\varOmega,\F,P)$ with a sigma--algebra $\F$, and a filtration $\{\emptyset,\varOmega\}=\F_0 \subset \F_1 \subset\dots\subset\F_T = \F$, we consider a stochastic process $\{W_t,~ t=0,\dots,T\}$ adapted to $\{\F_t,~t=0,\dots,T\}$.  Furthermore, we adopt the convention that any variable indexed by $t$ is $\F_t$--measurable. Additionally, we denote the sets of possible realizations of $W_t$ with $\Omega_t ,~ t=1,\dots, T$. The sets $\Omega_t ,~ t=1,\dots, T$ correspond to nested partitions of $\Omega$ given by the filtration $\{\F_t,~t=1,\dots,T\}$. We assume that each sample set $\Omega_t$ has a finite number of elements that is small enough to be enumerated computationally.\\

The \emph{post-decision information state} $I_t^x$ represents all the information in $S_t^x$ that is not in $R_t^x$. $I_t^x$ depends on the underlying exogenous random process, and contains only the information necessary to model the random transition from the current realization $\omega_t$ at time $t$ to the next random realization $\omega_{t+1}$ at time $t+1$. Hence, in the case of a Markov (lag 1) model, $I_t^x$ is given by a probability distribution
\bns
I^x_t = \Pb(\omega_{t+1}|S_t) = \Pb_{t+1}(\omega_{t+1}|\omega_t).
\ens
\begin{remark}
Please note that if the given problem features stagewise independence, i.e. a memoryless stochastic process, then all possible realizations $\omega_t$ would share the same
post-decision information state $I^x_t$.
\end{remark}

\textbf{Transition Function}\\
The amount of energy in the storage devices $R^B_{t}$ is adjusted to account for injected and withdrawn power,
\bn
R^B_{t+1, b}  &=& R^{B,x}_{t,b}\\
              &=& R^B_{t, b} + \eta^-_{b} x^{B-}_{t,b} - \eta_t^+ x^{B+}_{t,b}.
\en
The amount of power $E^W_t$ generated by the set of wind farms is adjusted as
\bns
E^W_{t+1} = E^W_{t} + \Ehat^W_{t+1}.
\ens
%where $\omega_{t+1}$ is the realization of the random vector containing the wind level updates.
% For each storage device $s$, the transition function for the energy in storage is given by:% %
% $$R_{t+1,s} = R_{t,s} + \eta_i^-x_{t}^-  - \frac{1}{\eta_i^+}x_t^+ $$
% For each wind power generator $g$, the transition function for the amount of power available:
% $$R_{t+1,g} = b_{t,g}$$
% In the case of stagewise independence, the information state
% The information state transition depends on the underlying assumptions: history dependence. How should we proceed?
% In the case of stagewise independence, the information state $I_t$ is empty, and no transition occurs.
% In the case of Markov uncertainty, the information state $I_{t+1}$ is a probability distribution on the set of realizations at time $t+2$, $$ I_{t+1} =  \Pb_{t+2}(\cdot|\omega_{t+1}).$$
% In the case of Markov uncertainty, $I_t$ coincides with the post-decision information state $I_t^x$, i.e. it is all the information
% The change in the information state reflects change in the amount of renewable energy which is adjusted by a random amount $\hat{I}_{t+1}$, $$I_{t+1} = I_t + \hat{I}_{t+1}.$$

The \emph{post-decision resource state} $R_t^x$ consists of two subvectors $R^{B,x}_{t}$ and $R^{G,x}_{t}$. The post-decision amount of energy in the batteries  $R^{B,x}_{t}$ is equal to the pre-decision amount $R^B_{t}$ adjusted for charging, discharging, as well as battery efficiencies,

Moreover, the post-decision state of the generators equals the power generation levels $x_{t}^{G+}$,
\bns
R^{G,x}_{t,g} = x^{G+}_{t,g},
\ens
from which we obtain the next pre-decision state
\bns
R^{G}_{t+1,g} = R^{G,x}_{t,g}.
\ens

\textbf{Objective Function}\\
We define the generation costs as the linear functions,
\bns
C(S_t, x_t) &=&  \langle c_t, x_t\rangle \\
            &=& \sum_{g\in\Gcal} Z^G_{t,g} c^G_{t,g} x^{G}_{t,g} + \sum_{b\in\Bcal} c^B_{t,b}x^{B}_{t,b}.
\ens
Our goal is to compute an optimal policy $X^{\pi^*}_t(S_t)$ that minimizes the total expected generation cost aggregated over the entire time horizon
\bn
\min_{\pi\in\Pi}\E\Big[\sum_{t=0}^T C(S_t, X_t^\pi(S_t)) | S_0 \Big]. \label{eq:objective}
\en
The expectation is taken with respect to a probability measure describing the possible outcomes of the energy from wind.
We note that in our energy application, the problem is highly time dependent, which is the reason we have indexed the policy itself by $t$, rather than just make it dependent on the state $S_t$ which depends on time.  Our experimental work will focus on modeling a problem over a full daily cycle in 5-minute increments, producing a problem with 288 time periods.

% The state of the system at time $t$ can be represented as $S_t = (R_t, I_t)$ where
% \begin{itemize}
% \item $R_t$ is the resource state, which is a vector of the amounts of energy available in all storage devices at time $t$ in MWh.
% \item $I_t$ denotes the information state, also known as state of the world. That is all the information that is exogenous to our model but still affects the
% outcome of our decisions.
% For our purposes, $I_t$ represents the wind energy available at all wind farms at time $t$ in MWh. We use $\Ic_t$ to denote the set of all possible information states $I_t$,
% and assume that $\Ic_t$ is discrete and finite.
% % We denote with $\Rc_t$ the set of all possible resource states $R_t$. is a continuous and bounded subset $\Rc_t\subset\R^m$.
% \end{itemize}
% We use $x_{tij}$ to denote the amount of energy transferred from $i$ to $j$ at time $t$.
%
% We would use subscript $g$ to stand for grid, and $d$ for demand.
% Further, we denote the set of wind farms with $\Jc^W$ and the set of storage devices with $\Jc^R$.
% Each storage device $j\in \Jc^R$ is characterized by the following specifications
% \begin{itemize}
% \item Manage multiple storage devices in a coordinated way.
% \item A capacity $R^c_j$ which is the maximum amount of energy the device can store in MWh.
% \item Maximum charging and discharging rates.
% \item Charging and discharging efficiency.
% \end{itemize}

\section{Stochastic optimization methods}
\label{sec:gridpaperadp}
We now contrast two algorithmic strategies for solving these problems which both approach the problem of approximating the value functions in Bellman's equation.  The first, widely known as the Stochastic Dual Dynamic Programming (SDDP) within the stochastic programming community was originally introduced in \cite{pereira1991multi}.
It uses multidimensional Benders cuts to approximate the value of being in a resource state $R_t$ using a {\it sampled} uncertainty model.  The second method uses the
language and notation of approximate dynamic programming, and approximates the value function using separable, piecewise linear approximations using a full uncertainty model.

These algorithms share the same fundamental style; as we show below, they are both a form of approximate dynamic programming, distinguished primarily by how they approximate the value function, and how they represent the underlying probability space.  However, there are several structural differences which make for an interesting comparison, especially in the setting of grid level storage, where the number of storage devices could range from single digits to hundreds, if we wish to anticipate a world of high penetration of renewables and lower cost batteries (which might be in the form of aggregators for electric vehicles).

\subsection{From stochastic programming to dynamic programming}
We begin by comparing the canonical models for SDDP and ADP.  SDDP, with its roots in the stochastic programming community, is a method for solving a problem that is often written in the form
\bn
% \min_{\substack{x_0 \in \Xcal_0}} \langle c_0, x_0\rangle + \E\left[\min_{\substack{x_0 \in \Xcal_1}} \langle c_1, x_1\rangle + \E\left[\dots +
% \E\left[\min_{\substack{x_T \in \Xcal_T}} \langle c_{T}, x_T\rangle \F_{T-1},x_{T-1}\right]\dots \F_1,x_1\right] |\F_0,x_0\right]. \nonumber \\
%  & & \label{eq:stochprogbasemodel}
\min_{\substack{x_0 \in \Xcal_0}(S_0)} \langle c_0, x_0\rangle + \E_1\left[\min_{\substack{x_1 \in \Xcal_1}(S_1)} \langle c_1, x_1\rangle + \E_2\left[\dots +
\E_T\left[\min_{\substack{x_T \in \Xcal_T}(S_T)} \langle c_{T}, x_T\rangle \right]\dots \right] \right]. \nonumber \\
 & & \label{eq:stochprogbasemodel}
\en
where $\Xcal_0(S_0) = \{x_0:A_0x_0 = b_0\}$, and for $t\geq 1$, we define the feasible sets as $\Xcal_t(S_t) = \{x_t: A_tx_t = b_t - B_{t-1}x_{t-1}, x_t\geq0\}$.
Here, it is assumed that $A_t$, $B_t$ and $b_t$, as well as the cost vector $c_t$, evolve randomly over time and describe the information contained in the stochastic process $\{W_t\}_{t=1}^T$, i.e. $(A_t, B_t, b_t, c_t)$ are $\F_t$-measurable matrices and vectors.
The stochastic programming community uses several notational styles, but one popular system
defines $\xi_t$ as the new information $(A_t, B_t, b_t, c_t)$, and lets $\xi_{[t]}$ be the history $\xi_1, \ldots, \xi_t$ (see \cite{ShDeRu2014}).
In this work, we use $\omega_t$ instead of $\xi_t$, and we denote our history by $\omega_{[t]}$.

It is common in the stochastic programming literature to write the state variable at time $t$ as $(x_{t-1},\omega_{[t]})$ which captures the dependence of the resource state $R_t$ on the decision $x_{t-1}$ (additional random inputs may be contained in $\omega_t$).  It is then possible to write a Bellman-style recursion as
\bn
Q_t(x_{t-1},\omega_{[t]}) = \min_{x_t} \big(c_t x_t + \E\big[Q_{t+1}(x_t,\omega_{[t+1]})|\omega_{[t]}\big]\big). \label{eq:stochprogbellman}
\en
The expectation is computationally intractable, but it is possible to replace it with a series of cuts (\cite{HiSe91}, \cite{pereira1991multi}), producing the linear program
\bn
Q_t(x_{t-1},\omega_{[t]}) = \min_{x_t\in\Xcal_t(x_{t-1}, \omega_{[t]}),v} (c_t x_t + v), \label{eq:stochprogbellman2}
\en
where
\bn
v \geq \alpha^k_{t+1}(\omega_{[t]}) + \beta^k_{t+1}(\omega_{[t]}) x_t, ~~~\mbox{for}~ k=1,\ldots, K, \label{eq:stochprogbellman3}
\en
and where $\Xcal_t$ captures the feasible region for $x_t$.  Here, equation \eqref{eq:stochprogbellman3} is generated by solving the dual problem for time $t+1$, which means that $K$ depends on the number of iterations that have been executed.  The indexing of the cuts in \eqref{eq:stochprogbellman3} reflects the fact that we are approximating the value at time $t+1$, but it is more accurate to say that it is approximating the recourse function around the post-decision state $S^x_t$ at time $t$ (the indexing of time should always reflect the information content of a variable when modeling a stochastic system).  

There are two computational issues with using the notational system of writing the state as $(x_{t-1},\omega_{[t]})$.  First, $x_{t-1}$ is generally a very high dimensional vector.  In the setting of our energy storage problem, it would have approximately 10,000 variables, one for each transmission line in the PJM grid.  Second, indexing on the history $\omega_{[t]}$ is, of course, problematic.  Even if our random variables were scalars (for each time period), the shortest horizon that we are going to consider in our work is 288 time periods (5-minute increments over 24 hours).  Retaining a history with more than three or four time periods (even with one variable per time period) is computationally intractable.

Using our notation, we would write
\bns
\E\big[Q_{t+1}(x_t,\omega_{[t+1]})|\omega_{[t]}\big] = V^x_t(S^x_t) = V^x_t(R^x_t,I^x_t).
\ens
This allows us to write our Bellman equation in the form of
\bns
V_t(S_t) = \min_{x_t\in\Xcal_t(S_t)} \big(c_t x_t + V^x_t(S^x_t)\big)
\ens
where $S_t = (R_t,I_t)$ is the pre-decision state, and $S^x_t = (R^x_t,I^x_t)$ is the post-decision state. Since our problem has a deterministic resource transition process, we know that $R_{t+1} = R^x_t = B_t x_t$.  We note that the dimension of $R_t$ and $R^x_t$ is equal to the number of storage devices, which might be as small as 1 or 10, or as large as 100 (in the experiments that we run).  Even a 10-dimensional resource vector would be too large if we were using a lookup table representation, but this is where convexity and Benders cuts allows us to obtain accurate approximations without enumerating the state.

This leaves the problem of the information state.  SDDP assumes that the process $W_t$ is a zero-th order Markov process (widely referred to as intertemporal or interstage independence), which means that $W_{t+1}$ is independent of $W_{t}$.  Under this assumption (which might at least be a reasonable approximation), the post-decision information state $I^x_t$ is empty and can be ignored, which means that $V^x_t(S^x_t) = V^x_t(R^x_t)$.  This represents a dramatic simplification of what is otherwise a very high-dimensional ($x_t$ has 10,000 dimensions) stochastic optimization problem with a long horizon (at least 288 time periods).

A form of the objective function \eqref{eq:stochprogbasemodel} that is more familiar to the stochastic programming community is to create a sampled ${\hat \Omega}\in\Omega$ and write
\bns
\min_{x_0, \ldots, x_T} &\sum_{\omega\in{\hat \Omega}}& \sum_{t=0}^T c_t(\omega)x_t(\omega) \label{eq:stochprogbasemodel2}\\
\mbox{s.t.~}A_0x_0 &=& b_0\\
% \ens
% and for $t=1,\ldots,T$, $\omega\in{\hat \Omega}$:
% \bns
B_{t-1}(\omega)x_{t-1}(\omega) + A_t(\omega) x_t(\omega) &=& b_t(\omega),~ t=1,\ldots,T,~ \omega\in{\hat \Omega}\\
x_t(\omega) &\geq& 0,~t=1,\ldots,T,~ \omega\in{\hat \Omega}
\ens
The notation $x_t(\omega)$ makes it possible for a decision at time $t$ to have access to the entire sample path $\omega$.  For this reason, we have to introduce {\it nonanticipativity constraints}.  This can be done by defining a history $h_t = (W_1, W_2, \ldots, W_t)$, with the set of histories $\Hcal_t = \{h_t(\omega), \omega\in{\hat \Omega}\}$.  Next let
\bns
{\hat \Omega}_t(h_t) = \{\omega\in{\hat \Omega}: (W_1(\omega),\ldots,W_t(\omega)) = h_t\},
\ens
be the set of all sample paths sharing the history $h_t$.  The nonanticipativity constraints are now given by
\bn
x_t(h_t) = x_t(\omega),~~\forall \omega\in{\hat \Omega}_t(h_t),~~h_t\in\Hcal_t. \label{eq:nonanticipativity}
\en
This is an example of a {\it sampled model} which is popular in stochastic optimization.  The standard approach is to generate ${\hat \Omega}$ in the form of a scenario tree, where histories are constructed by starting with a history $h_t$ at time $t$, and then branching by sampling realizations of $W_{t+1}$ given $h_t$.

We can make the transition to the formulation used in the dynamic programming community by replacing $x_t(\omega)$ with $x_t(h_t)$, which avoids the need for the nonanticipativity constraints \eqref{eq:nonanticipativity}.  Stochastic programmers will recognize $h_t$ as a node in the scenario tree.  We take this one step further by recognizing that we generally do not need the entire history.  Instead, we use the state $S_t$ which is the minimally dimensioned function of history $h_t$ which is necessary and sufficient (along with the exogenous information) to model our system from time $t$ onward.  We then have to choose $x_t(S_t)$ that satisfies the constraints $A_tx_t=b_t - B_{t-1}x_{t-1} = R_t$.

The dynamic programming community approaches the problem by recognizing that $x_t(S_t)$ is a function called a policy, that we write $X^\pi_t(S_t)$.  The objective function can now be written
\bn
\min_\pi \E \left[\sum_{t=0}^T C(S_t,X^\pi_t(S_t)) | S_0 \right].\label{eq:DPbasemodel}
\en
The dynamics of the system are captured through the transition function
\bn
S_{t+1} = S^M(S_t,x_t,W_{t+1}) \label{eq:transition}
\en
where $x_t = X^\pi_t(S_t)$.  We note that if we fix a policy, we can simulate the value of the policy for $\omega\in\Omega$ using
\bn
v^\pi(\omega) = \sum_{t=0}^T C(S_t(\omega),X^\pi_t(S_t(\omega))) \label{eq:simulatepolicy}
\en
where $S_{t+1}(\omega) = S^M(S_t(\omega),X^\pi_t(S_t(\omega)),W_{t+1}(\omega))$.  We note that while simulating the policy in \eqref{eq:simulatepolicy}, it is quite easy to make the outcome $W_{t+1}(\omega)$ dependent on both $S_t$ and $x_t$, which can be a valuable feature in energy applications (discharging energy into the grid can dampen electricity prices).  The decisions $x_t = X^\pi_t(S_t)$ must satisfy $x_t\in\Xcal_t = \{x_t: A_tx_t = b_t - B_{t-1}x_{t-1}, x_t\geq0\}$.  Nonanticipativity is satisfied automatically by writing the policy as dependent on the state rather than the sample path.

The dynamic programming community often skips equation \eqref{eq:DPbasemodel} and goes directly to Bellman's equation, which would be written
\bn
V_t(S_t) = \min_{x\in\Xcal_t} \big(C(S_t,x) + \E \big[ V_{t+1}(S_{t+1})|S_t\big]\big).  \label{eq:stochprogbellmanV}
\en
Exploiting the post-decision state allows us to drop the expectation, giving us
\bn
V_t(S_t) = \min_{x\in\Xcal_t} \big(C(S_t,x) + V^x_t(S^x_t)\big).  \label{eq:stochprogbellmanpost}
\en
If we make the same assumption that the information process $\{W_t\}_{t=1}^T$ is independent over time, then $S^x_t = R^x_t = B_tx_t$.  Exploiting the convexity of $V^x_t(R^x_t)$ (as well as $V_t(R_t,I_t)$), we can approximate $V^x_t(S^x_t)=V^x_t(R^x_t)$ using Benders cuts (as we did in equations \eqref{eq:stochprogbellman2}-\eqref{eq:stochprogbellman3}), but we can also experiment with other approximations.  Below, we test the idea of approximating the value function using separable, piecewise linear functions which we can write
\bn
V^x_t(R^x_t) \approx \sum_{i\in\Ical} \Vbar^x_{ti}(R^x_{ti}), \label{eq:separable}
\en
where $\Vbar^x_{ti}(R^x_{ti})$ is a one-dimensional piecewise linear and convex function.  These approximations can be estimated from dual variables from the sampled problem solved at time $t+1$ using algorithms such as CAVE and SPAR (see \cite{Po11}[Chapter 13], or \cite{Godfrey2001a} and \cite{PoRuTo04}).

We note that while we can assume that the post-decision information state $I^x_t$ is empty, there may be applications where the future depends on a compact information state.  For our energy application, we might characterize the weather using a small number of states that capture temperature and the likelihood of precipitation.  If we can represent these ``states of weather'' using perhaps 10 or 20 values, then we can create indexed convex approximations using either Benders cuts or
separable approximations.  Such models have been described as ``Markov'' in the literature (for some reason, many authors assume that a ``Markov'' model state space has to be small and discrete).  \cite{LoWoMi2015} and \cite{Lohndorf2015} pursue this idea under the name ``Approximate dual dynamic programming,'' but an alternative name is "Markov SDDP" which would help to communicate the method to the stochastic programming community.
%We use the term ``Markov SDDP'' which we think is more descriptive.  Also, all of these methods are basically a form of approximate dynamic programming, since they use Bellman's principle and develop approximations of the value function.

We are now going to present and test two algorithms that are structurally very similar. 
Recognizing that the fundamental structure of these two algorithms is quite similar, we undertake a series of qualitative and empirical comparisons that highlight the two key differences: how we model the information process (sampled or full), and how we approximate the value function (multidimensional Benders cuts or piecewise linear, separable).
% The first, SDDP, was developed within the stochastic programming community.  The second, which we refer to as ADP-SPWL (separable, piecewise linear), was developed using the language and notation of approximate dynamic programming.  There are important differences between the algorithms, but their overall structure is the same.

\subsection{SDDP with Benders cuts}
Stochastic dual dynamic programming (SDDP) has long enjoyed special attention from the stochastic programming community. The algorithm is summarized
% \ref{a:reg-sddp-markov}
in the online supplement.  Key features of the algorithm include:

\begin{itemize}
\item The algorithm solves a sampled version of the problem, using ${\hat \Omega}_t$ that are chosen before the algorithm starts.
\item The value functions are approximated by multidimensional Benders cuts.
\item It uses a forward pass, simulating the process of making decisions by solving a sequence of linear programs given by \eqref{eq:stochprogbellman2}-\eqref{eq:stochprogbellman3}.  Regularization is used to stabilize the solution.
\item Dual solutions are computed for {\it each} sample realization in ${\hat \Omega}_t$, computed in a backward pass. A new cut is created by averaging across the duals using the sample ${\hat \Omega}_t$.
\end{itemize}

\subsection{ADP-SPWL with separable, piecewise linear value function approximation}
The algorithm ADP-SPWL is described in the online supplement.
% \ref{a:adp}.
Key characteristics of the algorithm include:
\begin{itemize}
\item The value functions are approximated using separable, piecewise linear value functions.
\item It uses a forward pass, simulating the process of making decisions by solving a sequence of linear programs given by \eqref{eq:stochprogbellmanpost}-\eqref{eq:separable}. Samples at time $t$ are generated on the fly from the full sample space $\Omega_t$, which means that they can reflect the state $S_t$.
\item Numerical derivatives (requiring the solution of a linear program) are computed for {\it each} storage device.  These derivatives are then smoothed to create updated VFAs for each storage device.
\end{itemize}
We note that our algorithm could have been implemented as a pure forward pass procedure, as has been done in transportation applications (see e.g., \cite{Topaloglu2006}).  In these applications, computing numerical derivatives for each storage device would not be necessary.  However, a pure forward pass algorithm can exhibit slow convergence when decisions at one point in time have an impact many time periods in the future, which is the case in our energy storage application.  Our battery storage example models a day (or more) in 5-minute increments.  We may, for example, wish to store energy at 3 pm to use at 9 pm, which is 72 time periods in the future. For this type of problem, a backward pass dramatically accelerates the learning over time.  However, this means that we need to know the flow augmenting path from an increment of energy in a battery at time $t$ in terms of how it impacts the resource state $R_{t+1}$.  This is the reason that numerical derivatives are necessary.

\subsection{A comparison}

\begin{center}\label{t:comparison}
\begin{tabular}{|l|l|}
\hline
Stochastic dual dynamic programming                                  & Separable piece-wise linear value functions\\
dynamic programming                              & value functions \\
\hline
\hline
A fixed sample is generated once and used        & New samples are drawn each iteration from \\
for all iterations (the sampled model).           & the original information model. \\
\hline
Solves a subproblem for each random           & Solves a subproblem for each post-decision  \\
realization $\omega_t \in \hat{\Omega}_t$ at  each time   & resource dimension $R_{t,m}^x$ at  each time \\
$t=0,\ldots, T$.                       & $t=0,\ldots, T$. \\
\hline
Multidimensional Benders cuts for VFAs.           & Separable, piecewise linear VFAs. \\
\hline
Growing set of hyperplanes.                       & Growing set of kinks in each VFA. \\
\hline
Lower bounds for sampled model.                   & No lower bounds.   \\
\hline
Stochastic process must be independent of        & Sampled outcomes may depend on the state. \\
the state.                                        &  \\
\hline
\end{tabular}
\end{center}

Table \ref{t:comparison} shows a side-by-side comparison of characteristics of SDDP and the ADP algorithm using a separable, piecewise linear value function. There are several differences which should be highlighted.
\begin{itemize}
\item SDDP requires the use of a sampled model because it is averaging multidimensional cuts.  The likelihood of visiting the same multidimensional resource state $R_t$ on two successive iterations is nearly zero.  By contrast, ADP-SPWL exploits the use of separable approximations, which allows us to smooth observations of slopes from different sample realizations into a single function. ADP-SPWL updates the marginal value functions for all of the resource dimensions at each time period, in each backward pass.
\item By generating multidimensional cuts, SDDP enjoys the feature that it can generate upper bounds on the solution, which can be used to help evaluate the quality of the solution.  However, for our energy storage application, it is important to use care when interpreting these bounds, because the objective function features a large constant term representing the value of a myopic policy (setting the value functions equal to zero).  When evaluating solution quality using upper and lower bounds, the gap can seem small if we do not subtract this constant term.
\item SDDP requires solving a linear program for each sample realization at each time period.  ADP-SPWL requires solving a linear program for each dimension of our resource state variable.
\item Since SDDP is solving a sampled approximation, the bounds that it generates are, strictly speaking, bounds on the optimal solution of the sampled problem. 
\item SDDP must use a sampled model that is generated before the algorithm starts.  ADP-SPWL draws samples dynamically during the forward pass, making it possible to generate samples that depend on the state $S_t$ and/or the action $x_t$.
\item Benders decomposition is known to show slow convergence as the dimensionality of the resource state grows, but the impact of dimensionality will depend on the characteristics of the problem.  The separable approximations for ADP-SPWL has been shown to work well at high dimensions, but this work is based on experiments in transportation and logistics where separability is likely to be a better approximations (see \cite{Topaloglu2006} and \cite{BeChPo2014} for illustrations).  In our grid application, substitution of energy between storage devices is much easier, suggesting that a separable approximation may not work as well.
\item Neither algorithm has been tested against serious competition. SDDP has been evaluated purely on the basis of its bounds on the sampled approximation.  ADP-SPWL has been evaluated primarily through comparisons against optimal solutions on deterministic problems (see \cite{ToPo06}) or scalar, stochastic problems \cite{JiPhPo14}.
\end{itemize}

\section{Computational comparisons of SDDP and ADP-SPWL}
\label{sec:gridpaperexperiments}

In this section we study the computational performance of the algorithms proposed above.
% We examine the performance of the methods for different dimensions of the value function approximations, as well as different sequences of regularization coefficients.
We introduce the following simplifications to the state space:
\begin{itemize}
\item The post-decision resource space is collapsed around $R^B_{t}$. In this way, we ignore rampup limits for the sake of numerical testing.
\item The information space is collapsed around $E^W_t$ only since that is the main source of volatility.
\end{itemize}
In our numerical experiments we focus our analysis on the following questions:
\begin{itemize}
% \item How is the computational performance of Algorithms \ref{a:reg-sddp-markov} and \ref{a:adp} affected by:
\item How is the computational performance of SDDP and ADP-SPWL affected by:
% \item We examine the impact on the computational performance of Algorithms \ref{a:reg-sddp} and \ref{a:reg-sddp-markov} caused by the following:
\begin{itemize}
\item the dimension of the resource vector $R_t^x$?
\item the size of the sample sets $|\hat{\Omega}_t|$?
\end{itemize}
% \item How does the performance of  Algorithm \ref{a:reg-sddp-markov} compare to its regularized counterpart?
\end{itemize}

Our experimental work was conducted using the setting of optimizing grid level storage for a large transmission grid managed by PJM Interconnection.  PJM manages grid level storage devices from a single location, making it a natural setting for testing our algorithms.  As of this writing, grid level storage is dropping in price, providing a meaningful setting to evaluate the performance of our algorithms for a wide range of storage devices, challenging the ability of the algorithms to
handle high dimensional applications.  For this reason, we conducted tests on networks with up to 100 storage devices.  These are much higher dimensional problems than prior research that has focused on the management of water reservoirs.
In order to be able to use high-dimensional Benders approximations, we consider a quadratic regularization extension of SDDP that was originally introduced by \cite{asamov2015regularized}. 
The algorithm was implemented in Java, and the IBM ILOG CPLEX 12.4  solver was used for the solution of both linear and quadratic convex optimization problems.
In addition, the relative complementarity tolerance of CPLEX was set to $10^{-12}$.

Another distinguishing feature of our grid storage setting (compared to prior experimental work) is that a natural time step is 5 minutes, which is the frequency with which real--time electricity prices (known as LMPs, for locational marginal prices) are updated on the PJM grid.  We anticipate using storage devices to hold energy over horizons of several hours.  For this reason, we used a 24 hour model, divided into 5--minute increments, for 288 time periods, which is quite large compared to many applications using this algorithmic technology.

Below we describe the construction of the network, the representation of the exogenous stochastic process, and finally we present the results of an extensive set of experiments investigating the effect of regularization, the number of storage devices (which determines the dimensionality of $R_t^x$), and the presence of an exogenous post-decision information state, on the rate of convergence and solution quality.

\subsection{The network}

We performed our experiments using an aggregated version of the PJM grid. Instead of the full network with 9,000 buses and 14,000 transmission lines, we limited our analysis to the higher voltage lines, producing a grid with 1,360 buses and 1,715 transmission lines.
% The power generators include 396 gas turbines (23,309 MW), 50 combined cycle generators (21,248 MW),
% 264 steam generators (73,374 MW), 31 nuclear reactors (31,086 MW), and 84 conventional hydro power generators (2,217 MW).
Off--shore wind power was simulated for a set of hypothetical wind turbines with a combined maximum capacity of 16 GW.  Moreover, we consider a daily time horizon with 5--minute discretization resulting in a total of 288 time periods.

The data was prepared by first running a unit--commitment simulator called SMART--ISO that determines which generators are on or off at each point in time, given forecasts of wind generated from a planned set of off--shore wind farms.  We made the assumption that the use of grid level storage would not change which generators are on or off at any point in time.  However, we simultaneously optimize the generator output levels, while charging and discharging of storage devices around the grid in the presence of stochastic injections from the wind farms.

We placed the distributed storage devices at the points--of--interconnection for wind farms, as well as the buses with the highest demand. Each storage device is characterized by its minimum and maximum energy capacity, its charging and discharging efficiency, and its variable storage cost. The control of multiple storage devices in a distributed energy system is a challenging task that depends on a variety of factors such as the location of each device, and the presence of transmission line
congestion. A good storage algorithm needs to respond to daily variations in supply, demand and congestion, taking advantage of opportunities to store energy near generation points (to avoid congestion) or near load points (during off--peak periods). It has to balance when and where to store and discharge in a stochastic, time--dependent setting, providing a challenging test environment for our algorithm.

\subsection{The exogenous information}
\label{s:7.1}

Our only source of uncertainty (the exogenous information) was from the injected wind from the offshore wind farms.  In order to calibrate our stochastic wind error model, we employed historical wind data and speed measurements of off--shore wind for the month of January 2010.  For each time period $t$, we consider a set $\hat{\Omega}_t$ of vectors of possible wind speed realizations which correspond to $|\hat{\Omega}_t|$ different weather regimes. Plots of simulated wind power at a given
wind farm can be seen in Figure 2
% \ref{fig:wind}.
in the online supplement.
% \input{wind.tex}

% In general, the exogenous information process can be characterized by one of the following: stagewise independence, compact state variables (so-called Markov processes), or scenario-dependence (path dependence). For some instances, the latter case could be reduced to one of the former two by applying an appropriate transformation. %as described in section \ref{s:5}.
% In our experiments, we consider instances with stagewise-independent transitions between ten equally likely scenarios.  When we assumed stagewise independence, we would sample from each of these 10 scenarios with equal probability at each time period. 
% For the problems with Markov uncertainty, we assumed that at every time period $t$, the probability of continuing with the same weather regime at time $t+1$ is 91 percent.  Additionally, each of the remaining nine regimes can be visited at time $t+1$ with a probability of 1 percent.

\subsection{Deterministic Experiments}
In the online supplement we present deterministic experiments that allow us to perform a side by side comparisons of objective values and decisions made by SDDP and ADP-SPWL. % , which provides an important test of both algorithms.
We can see that in the instance with only five storage devices, the SDDP solution is closer to the optimal than the solution of ADP-SPWL. However, as the dimension of the post-decision value functions increases, the separable value function approximations outperform their Benders counterparts in terms of both objective value and solution quality.
% This raises the following issue. If the SDDP approximations do not perform well in a deterministic setting, then how could we count on them to solve even more complex stochastic optimization problems? We address this question next.  
Still, the real test has to be on a stochastic dataset, since a deterministic problem allows the algorithms to learn the states of other storage devices at each point in time. 

\subsection{Stochastic Experiments}
 In this section, we turn our attention to stochastic experiments, and at the same time address issues related to the use of a sampled model by SDDP.
% Further, we performed parameter tuning as described in section \ref{s:tuning}. We set the relative complementarity tolerance of CPLEX to $10^{-12}$, and used a geometric regularization sequence with $\varrho^0 = 1$ and $r=0.95$. Moreover, the scaling matrices $Q_t, t=0,\dots,T$ are set to the identity matrix which implies that the amount of energy in each storage device has the same weight in the regularization term. In this section we examine the performance of Algorithms \ref{a:reg-sddp} and \ref{a:reg-sddp-markov} when the number of storage devices (dimension of the resource state variable) is $|R^x_t| = 50, 100, 200, 500$.
% In this section we study the computational behavior of the SDDP and ADP-SPLW. % Plots of the behavior of all methods can be found in Figures \ref{fig:25}, \ref{fig:50}, \ref{fig:100} below.  
% Each figure shows the results for Algorithm \ref{a:reg-sddp-markov} on the left, its regularized counterpart in the middle, and Algorithm \ref{a:adp} on the right.
The first question that we study is the computational effect of increasing the dimensionality of $R^x_t$.  This is a major issue, largely overlooked in the SDDP community.  SDDP has always been presented as a way of circumventing the curse of dimensionality, but in practice it has been used only for the solution of instances with very low dimensional value functions. Thus, the question of its practical applicability to high--dimensional problems has been left unanswered until
now.
The plots below illustrate several important points about solving such large-scale stochastic problems.
First, both the regularized version of SDDP, as well as ADP-SPWL seem applicable for instances with large resource dimensions $|R^x_t|$ that would be intractable for non-regularized Benders methods (even in the deterministic case).
% The use of regularization significantly accelarates the convergence of SDDP-type methods, and its usefulness grows as the number of dimensions increases.
The results suggest that SDDP regularization allows practitioners to consistently obtain high quality solutions within approximately 50 iterations for all of the given problems.
However, during initial iterations the ADP-SPWL approach can exhibit even faster convergence than the regularized SDDP.
% Second, we can see that the importance of the sample size $\hat{\Omega}_t$ is also diminished when the problem Benders-type approximations are regularized. Such a property can be very useful when the dimension of the information space increases beyond 4 or 5.
% In that case, practitioners might attempt to minimize the sampling error that would stem from using a sample $\hat{\Omega}_t$ instead of the full space $\Omega_t$.
% For high-dimensional random processes, that might entail sample sizes of thousands of realizations (or more), which could potentially present a major computational challenge since SDDP-type methods need to solve an optimization problem for each random realization at each time step in the backward pass of every iteration.
% However, our results indicate that the discrepancies stemming from using a small sample can be mitigated by regularization.
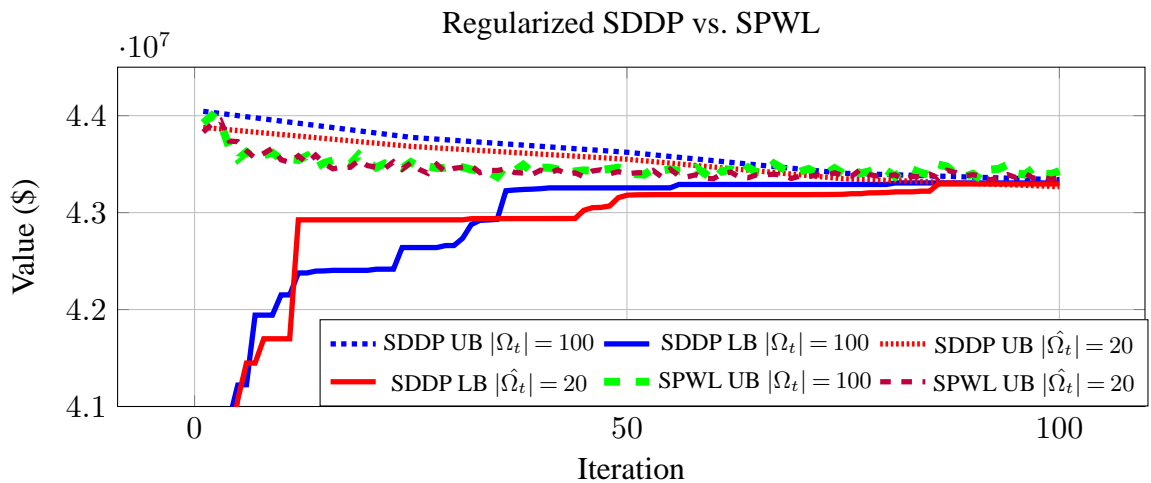
\begin{figure}[H]
\begin{tikzpicture}[baseline]
\pgfplotsset{every axis legend/.append style={font=\footnotesize, at={(0.6,0.0)},anchor=south, legend columns=3}}
\begin{axis}[
height=2.4in,
width=6.0in,
grid=major,
xtick = {0, 50,...,100},
ymin=41000000.0,
ymax=44500000.0,
title={Regularized SDDP vs. SPWL},
xlabel={Iteration},
ylabel={Value (\$)}]
\addplot[color=blue, dash pattern=on 2pt off 2pt, line width = 2] coordinates {
 (1, 44044957.5) (25, 43776422.7) (50, 43622435.4) (75, 43406009.0) (100, 43339521.6)
};
\addlegendentry{SDDP UB $|{\Omega_t}| = 100$}
\label{fig:25:m_sddp_ub}
\addplot[color=blue, line width = 2] coordinates {
 (1, 28530323.5) (2, 37789327.7) (3, 37789327.7) (4, 40862721.6) (5, 41221130.8) (6, 41224524.2) (7, 41942364.6) (8, 41942364.6) (9, 41942364.6) (10, 42150991.5) (11, 42150991.5) (12, 42377585.3) (13, 42377605.9) (14, 42396631.3) (15, 42398699.2) (16, 42404426.5) (17, 42404450.3) (18, 42404541.5) (19, 42404541.5) (20, 42404638.5) (21, 42416208.2) (22, 42416510.0) (23, 42416600.5) (24, 42639544.2) (25, 42639544.2) (26, 42639544.2) (27, 42639544.2) (28, 42639544.2) (29, 42659837.4) (30, 42660456.1) (31, 42735491.4) (32, 42877699.7) (33, 42920276.6) (34, 42924457.5) (35, 42930023.0) (36, 43225441.4) (37, 43233107.0) (38, 43238698.6) (39, 43241753.8) (40, 43246836.6) (41, 43255034.2) (42, 43255051.3) (43, 43255051.3) (44, 43255052.5) (45, 43255052.5) (46, 43255052.6) (47, 43255053.1) (48, 43255055.8) (49, 43255055.8) (50, 43255055.8) (51, 43255055.8) (52, 43255055.8) (53, 43255090.4) (54, 43255391.0) (55, 43255391.0) (56, 43290140.6) (57, 43290221.3) (58, 43290228.0) (59, 43290247.0) (60, 43290247.1) (61, 43290247.1) (62, 43290249.2) (63, 43290249.2) (64, 43290249.2) (65, 43290249.2) (66, 43290249.2) (67, 43290250.6) (68, 43290250.6) (69, 43290250.6) (70, 43290250.6) (71, 43290250.6) (72, 43290250.7) (73, 43290250.7) (74, 43290250.7) (75, 43290250.7) (76, 43290251.0) (77, 43290251.0) (78, 43290251.0) (79, 43290251.0) (80, 43290251.0) (81, 43305937.9) (82, 43305971.4) (83, 43305973.7) (84, 43305974.0) (85, 43305978.5) (86, 43305978.5) (87, 43305978.5) (88, 43305978.5) (89, 43305994.5) (90, 43305994.5) (91, 43305994.5) (92, 43305994.5) (93, 43305994.5) (94, 43305994.7) (95, 43305994.9) (96, 43305994.9) (97, 43305994.9) (98, 43305994.9) (99, 43305994.9) (100, 43305994.9)
};
\addlegendentry{SDDP LB $|{\Omega_t}|=100$}
\label{fig:25:m_sddp_lb}
\addplot[color=red, dash pattern=on 1pt off 1pt, line width = 2] coordinates {
 (1, 43879295.7) (25, 43683563.6) (50, 43550558.9) (75, 43349256.9) (100, 43268408.8)
};
\addlegendentry{SDDP UB $|\hat{\Omega_t}| = 20$}
\label{fig:25:m_reg_sddp_ub}
\addplot[color=red, line width = 2] coordinates {
 (1, 27663288.0) (2, 38943935.7) (3, 38962403.1) (4, 40500661.3) (5, 41072698.6) (6, 41448060.4) (7, 41448060.4) (8, 41698659.0) (9, 41698682.1) (10, 41698690.7) (11, 41698690.7) (12, 42925575.0) (13, 42925575.2) (14, 42925575.2) (15, 42925575.2) (16, 42925575.2) (17, 42925575.2) (18, 42925575.2) (19, 42925575.2) (20, 42925575.2) (21, 42925575.2) (22, 42925575.2) (23, 42925575.2) (24, 42925575.2) (25, 42925575.2) (26, 42925575.2) (27, 42925575.2) (28, 42925575.2) (29, 42925575.2) (30, 42925575.2) (31, 42925575.2) (32, 42937369.8) (33, 42937744.7) (34, 42937839.3) (35, 42938008.8) (36, 42938057.9) (37, 42938076.0) (38, 42938080.2) (39, 42938094.7) (40, 42938217.1) (41, 42938508.5) (42, 42938628.4) (43, 42938743.0) (44, 42938954.3) (45, 43023597.3) (46, 43051049.4) (47, 43053844.0) (48, 43068460.7) (49, 43154330.3) (50, 43179716.4) (51, 43183330.0) (52, 43184623.1) (53, 43185240.2) (54, 43185301.2) (55, 43185321.0) (56, 43185326.4) (57, 43185339.8) (58, 43185359.3) (59, 43185360.1) (60, 43185360.8) (61, 43185363.4) (62, 43185363.5) (63, 43185363.5) (64, 43185363.5) (65, 43185363.5) (66, 43185363.5) (67, 43185371.1) (68, 43185371.1) (69, 43185371.1) (70, 43185371.1) (71, 43185371.1) (72, 43186349.4) (73, 43186349.4) (74, 43187526.2) (75, 43189024.9) (76, 43195301.1) (77, 43195301.1) (78, 43203847.6) (79, 43204893.9) (80, 43207099.5) (81, 43213761.8) (82, 43214284.0) (83, 43214488.3) (84, 43220806.4) (85, 43220806.4) (86, 43295627.5) (87, 43295628.2) (88, 43295628.2) (89, 43295629.2) (90, 43295633.9) (91, 43295640.3) (92, 43295640.3) (93, 43295640.4) (94, 43295642.5) (95, 43295645.1) (96, 43295648.9) (97, 43295652.2) (98, 43295655.9) (99, 43295656.4) (100, 43295657.6)
};
\addlegendentry{SDDP LB $|\hat{\Omega_t}| = 20$}
\label{fig:25:m_reg_sddp_lb}
\addplot[color=green, dash pattern=on 6pt off 6pt, line width = 3] coordinates {
 (1, 43924227.3) (2, 44002091.9) (3, 43932288.9) (4, 43669717.9) (5, 43551556.2) (6, 43625306.0) (7, 43564273.0) (8, 43590679.6) (9, 43617977.1) (10, 43550408.5) (11, 43550273.4) (12, 43543178.0) (13, 43586206.3) (14, 43465866.5) (15, 43513326.2) (16, 43539279.6) (17, 43535768.0) (18, 43499264.8) (19, 43598764.0) (20, 43523224.9) (21, 43464171.3) (22, 43508129.6) (23, 43485193.4) (24, 43459376.2) (25, 43521375.0) (26, 43473411.2) (27, 43476901.4) (28, 43459023.0) (29, 43461956.7) (30, 43487611.1) (31, 43467909.0) (32, 43444541.8) (33, 43427907.6) (34, 43417941.5) (35, 43369067.8) (36, 43478195.6) (37, 43431617.8) (38, 43493931.7) (39, 43417961.7) (40, 43463596.0) (41, 43463388.6) (42, 43445588.7) (43, 43408755.1) (44, 43425995.3) (45, 43436004.5) (46, 43431935.6) (47, 43440716.2) (48, 43451403.6) (49, 43464193.0) (50, 43440225.9) (51, 43397368.5) (52, 43390865.4) (53, 43451394.9) (54, 43432748.1) (55, 43447172.0) (56, 43430718.0) (57, 43424445.1) (58, 43464186.9) (59, 43400871.1) (60, 43503518.7) (61, 43422230.4) (62, 43448871.3) (63, 43387458.0) (64, 43466364.6) (65, 43417752.5) (66, 43417529.2) (67, 43506815.7) (68, 43454123.2) (69, 43406367.1) (70, 43422096.4) (71, 43450449.5) (72, 43403076.3) (73, 43400233.8) (74, 43458227.3) (75, 43404401.3) (76, 43452885.8) (77, 43404940.7) (78, 43425219.4) (79, 43460168.6) (80, 43436375.8) (81, 43371105.4) (82, 43461461.3) (83, 43445584.9) (84, 43411362.0) (85, 43386915.7) (86, 43497266.9) (87, 43475061.8) (88, 43400480.6) (89, 43367922.4) (90, 43400581.9) (91, 43338236.4) (92, 43444221.2) (93, 43451974.3) (94, 43494440.9) (95, 43406412.5) (96, 43392541.1) (97, 43459780.0) (98, 43402200.8) (99, 43396948.8) (100, 43428206.1)
};
\addlegendentry{SPWL UB $|{\Omega_t}| = 100$}
\label{fig:25:reg_adp_ub}
\addplot[color=purple, dash pattern=on 4pt off 4pt, line width = 2] coordinates {
 (1, 43830901.9) (2, 43913172.0) (3, 43951009.6) (4, 43738606.2) (5, 43732007.0) (6, 43638757.5) (7, 43561367.7) (8, 43610890.0) (9, 43663194.0) (10, 43544771.9) (11, 43518257.0) (12, 43585384.4) (13, 43510951.0) (14, 43592557.1) (15, 43504918.0) (16, 43451845.0) (17, 43469240.3) (18, 43467380.3) (19, 43539830.7) (20, 43484359.5) (21, 43469128.9) (22, 43457144.9) (23, 43523921.5) (24, 43447227.5) (25, 43514483.7) (26, 43445575.8) (27, 43404278.8) (28, 43498141.3) (29, 43419761.4) (30, 43456626.2) (31, 43468070.2) (32, 43482014.7) (33, 43412405.1) (34, 43435079.2) (35, 43434543.0) (36, 43450030.3) (37, 43413993.1) (38, 43407869.3) (39, 43392773.8) (40, 43437490.2) (41, 43406660.2) (42, 43442525.3) (43, 43396805.9) (44, 43410221.8) (45, 43369680.5) (46, 43428865.0) (47, 43441861.1) (48, 43399724.3) (49, 43422837.8) (50, 43399567.3) (51, 43409850.7) (52, 43412528.3) (53, 43435532.2) (54, 43370133.7) (55, 43446912.8) (56, 43380983.9) (57, 43370272.2) (58, 43390805.5) (59, 43355118.4) (60, 43348987.5) (61, 43404044.4) (62, 43417538.1) (63, 43365899.7) (64, 43359408.8) (65, 43389861.3) (66, 43390172.9) (67, 43395943.1) (68, 43385529.0) (69, 43352220.4) (70, 43375685.3) (71, 43371285.4) (72, 43388061.2) (73, 43382154.9) (74, 43410307.6) (75, 43380067.7) (76, 43376584.3) (77, 43444603.5) (78, 43392168.4) (79, 43376934.5) (80, 43379406.5) (81, 43389379.0) (82, 43340904.7) (83, 43373316.8) (84, 43367929.9) (85, 43397023.1) (86, 43356211.5) (87, 43436658.4) (88, 43415309.7) (89, 43379481.4) (90, 43390459.9) (91, 43374948.2) (92, 43336765.3) (93, 43362529.6) (94, 43392089.2) (95, 43360680.8) (96, 43407010.8) (97, 43353488.3) (98, 43367991.2) (99, 43343962.8) (100, 43360742.5)
};
\addlegendentry{SPWL UB $|\hat{\Omega_t}| = 20$}
\label{fig:25:reg_adp_ub}
\end{axis}
\end{tikzpicture}
\caption{Numerical comparison of multistage stochastic optimization methods for $|R^x_t|= 25$.}
\label{fig:25}
\end{figure}
\begin{figure}[H]
\begin{tikzpicture}[baseline]
\pgfplotsset{every axis legend/.append style={font=\footnotesize, at={(0.6,0.0)},anchor=south, legend columns=3}}
\begin{axis}[
height=2.4in,
width=6.0in,
grid=major,
xtick = {0, 50,...,100},
ymin=41000000.0,
ymax=44500000.0,
title={Regularized SDDP vs. SPWL},
xlabel={Iteration},
ylabel={Value (\$)}]
\addplot[color=blue, dash pattern=on 2pt off 2pt, line width = 2] coordinates {
 (1, 43890303.2) (25, 43700397.7) (50, 43600444.5) (75, 43448789.1) (100, 43402142.9)
};
\addlegendentry{SDDP UB $|{\Omega_t}| = 100$}
\label{fig:50:m_sddp_ub}
\addplot[color=blue, line width = 2] coordinates {
 (1, 20373557.8) (2, 41750673.1) (3, 41752329.3) (4, 41753679.9) (5, 41755551.7) (6, 41757835.4) (7, 41759912.2) (8, 41762346.0) (9, 41763666.4) (10, 41765069.7) (11, 41766754.4) (12, 41768461.8) (13, 41772038.7) (14, 41779103.6) (15, 42239084.2) (16, 42271719.1) (17, 42276057.4) (18, 42292403.5) (19, 42303269.5) (20, 42321665.8) (21, 42347375.5) (22, 42350669.7) (23, 42353176.0) (24, 42353875.7) (25, 42365332.7) (26, 42365599.9) (27, 42365858.5) (28, 42367830.2) (29, 42368003.0) (30, 42368144.2) (31, 42368335.7) (32, 42368585.1) (33, 42368800.8) (34, 42368902.2) (35, 42372799.7) (36, 42373020.5) (37, 42373593.0) (38, 42373858.5) (39, 42374025.7) (40, 42374247.1) (41, 42374454.6) (42, 42374492.6) (43, 42374503.3) (44, 42374545.2) (45, 42573942.1) (46, 42651353.5) (47, 42652034.1) (48, 42696335.8) (49, 42734439.7) (50, 42734883.8) (51, 42738662.0) (52, 42740978.1) (53, 42742156.6) (54, 42750061.1) (55, 42752786.3) (56, 42763733.3) (57, 42763743.3) (58, 42763747.6) (59, 42763907.7) (60, 42763936.6) (61, 42763987.5) (62, 42764013.6) (63, 42764082.3) (64, 42764091.3) (65, 42764097.7) (66, 42764098.2) (67, 42764103.8) (68, 42764108.0) (69, 42764109.5) (70, 42764109.9) (71, 42764111.4) (72, 42764163.3) (73, 42764184.1) (74, 42764191.5) (75, 42764193.0) (76, 42764193.7) (77, 42764195.2) (78, 42764238.1) (79, 42764239.5) (80, 42764244.0) (81, 42764245.3) (82, 42764248.0) (83, 42764249.1) (84, 42764261.5) (85, 42764262.8) (86, 42764286.1) (87, 42764287.5) (88, 42764293.7) (89, 42764295.5) (90, 42764296.3) (91, 42764302.5) (92, 42764326.7) (93, 42764381.9) (94, 42764387.5) (95, 42764390.7) (96, 42764407.9) (97, 42764558.5) (98, 42764560.2) (99, 42765182.8) (100, 42765183.7)
};
\addlegendentry{SDDP LB $|{\Omega_t}|=100$}
\label{fig:50:m_sddp_lb}
\addplot[color=red, dash pattern=on 1pt off 1pt, line width = 2] coordinates {
 (1, 44055102.4) (25, 43665620.4) (50, 43554661.7) (75, 43434794.7) (100, 43436716.9)
};
\addlegendentry{SDDP UB $|\hat{\Omega_t}| = 20$}
\label{fig:50:m_reg_sddp_ub}
\addplot[color=red, line width = 2] coordinates {
 (1, 20339147.1) (2, 38669580.2) (3, 40346603.0) (4, 42248697.1) (5, 42250231.5) (6, 42267746.7) (7, 42271517.1) (8, 42275916.3) (9, 42280099.3) (10, 42284129.4) (11, 42288803.2) (12, 42293986.8) (13, 42298180.9) (14, 42302904.3) (15, 42306077.4) (16, 42307535.1) (17, 42312227.7) (18, 42316519.9) (19, 42322572.3) (20, 42337323.2) (21, 42342475.7) (22, 42358102.9) (23, 42360670.2) (24, 42362484.9) (25, 42363004.0) (26, 42363370.9) (27, 42370808.6) (28, 42372723.4) (29, 42373006.0) (30, 42373272.1) (31, 42373533.2) (32, 42373633.5) (33, 42375942.3) (34, 42382469.0) (35, 42386450.0) (36, 42417075.8) (37, 42419819.8) (38, 42424147.0) (39, 42424300.8) (40, 42527554.0) (41, 42532581.7) (42, 42582229.9) (43, 42670302.0) (44, 42704759.6) (45, 42750468.0) (46, 42752149.5) (47, 42752218.1) (48, 42752358.7) (49, 42761231.6) (50, 42762150.2) (51, 42766357.4) (52, 42769702.2) (53, 42770717.5) (54, 42770775.9) (55, 42770802.5) (56, 42770859.6) (57, 42770861.6) (58, 42770949.1) (59, 42770990.8) (60, 42771033.0) (61, 42772717.1) (62, 42772748.3) (63, 42773220.9) (64, 42773282.9) (65, 42773299.8) (66, 42773434.6) (67, 42773558.0) (68, 42773670.5) (69, 42773680.9) (70, 42774990.6) (71, 42775037.2) (72, 42775111.4) (73, 42775585.1) (74, 42775636.1) (75, 42776073.6) (76, 42776074.3) (77, 42776672.8) (78, 42776894.6) (79, 42776978.1) (80, 42777095.4) (81, 42777177.9) (82, 42777198.9) (83, 42777229.0) (84, 42777232.2) (85, 42777247.3) (86, 42777271.8) (87, 42777315.4) (88, 42777318.6) (89, 42777364.4) (90, 42777508.6) (91, 42777585.9) (92, 42777589.5) (93, 42777656.8) (94, 42777664.9) (95, 42777758.2) (96, 42777936.6) (97, 42777965.2) (98, 42778016.8) (99, 42778027.6) (100, 42778029.0)
};
\addlegendentry{SDDP LB $|\hat{\Omega_t}| = 20$}
\label{fig:50:m_reg_sddp_lb}
\addplot[color=green, dash pattern=on 6pt off 6pt, line width = 3] coordinates {
 (1, 43918045.6) (2, 43887234.8) (3, 43852088.2) (4, 43776626.2) (5, 43766899.9) (6, 43695631.0) (7, 43623589.3) (8, 43506089.2) (9, 43512753.4) (10, 43542627.9) (11, 43490327.5) (12, 43514227.1) (13, 43483887.2) (14, 43493984.3) (15, 43478264.1) (16, 43468857.7) (17, 43453771.4) (18, 43461477.2) (19, 43441773.5) (20, 43400619.0) (21, 43377515.2) (22, 43458544.5) (23, 43402127.9) (24, 43459635.9) (25, 43653450.9) (26, 43548664.6) (27, 43477308.4) (28, 43539740.4) (29, 43480068.6) (30, 43487184.8) (31, 43472835.6) (32, 43407702.8) (33, 43417834.0) (34, 43420152.7) (35, 43446650.8) (36, 43416384.0) (37, 43441419.5) (38, 43389753.2) (39, 43465907.6) (40, 43423820.8) (41, 43392328.9) (42, 43407816.9) (43, 43428471.6) (44, 43358725.1) (45, 43346749.9) (46, 43385131.7) (47, 43352492.7) (48, 43375782.7) (49, 43369532.3) (50, 43371114.8) (51, 43374083.7) (52, 43338986.5) (53, 43358041.2) (54, 43404739.9) (55, 43406242.2) (56, 43393637.0) (57, 43373306.1) (58, 43390873.1) (59, 43414675.3) (60, 43342031.2) (61, 43366924.2) (62, 43421401.2) (63, 43408337.0) (64, 43324744.9) (65, 43362268.5) (66, 43380084.2) (67, 43362068.3) (68, 43343675.6) (69, 43370918.8) (70, 43366535.8) (71, 43328210.7) (72, 43313959.9) (73, 43354222.8) (74, 43330128.5) (75, 43312305.1) (76, 43336423.5) (77, 43375979.3) (78, 43285618.2) (79, 43348134.5) (80, 43351752.0) (81, 43372192.1) (82, 43318146.6) (83, 43354466.4) (84, 43339536.3) (85, 43384898.4) (86, 43324460.0) (87, 43315310.7) (88, 43324186.6) (89, 43304250.5) (90, 43348030.9) (91, 43350884.3) (92, 43342419.7) (93, 43336588.0) (94, 43349590.0) (95, 43325635.6) (96, 43357536.1) (97, 43344787.5) (98, 43289942.9) (99, 43360970.3) (100, 43355966.1)
};
\addlegendentry{SPWL UB $|{\Omega_t}| = 100$}
\label{fig:50:reg_adp_ub}
\addplot[color=purple, dash pattern=on 4pt off 4pt, line width = 2] coordinates {
 (1, 43912056.7) (2, 43864678.2) (3, 43946927.6) (4, 43746174.7) (5, 43601985.0) (6, 43544863.0) (7, 43490023.4) (8, 43479687.8) (9, 43486452.0) (10, 43450783.1) (11, 43455504.1) (12, 43459967.3) (13, 43480555.8) (14, 43446436.7) (15, 43457278.7) (16, 43423755.2) (17, 43433162.0) (18, 43395650.6) (19, 43462474.1) (20, 43466092.7) (21, 43435563.7) (22, 43447701.6) (23, 43437123.9) (24, 43454070.0) (25, 43407178.4) (26, 43420655.7) (27, 43388966.4) (28, 43365931.0) (29, 43363557.5) (30, 43370145.6) (31, 43362505.0) (32, 43398446.8) (33, 43384668.4) (34, 43366387.3) (35, 43360846.2) (36, 43380336.2) (37, 43361117.8) (38, 43316801.1) (39, 43351348.0) (40, 43377178.2) (41, 43324812.4) (42, 43382142.1) (43, 43327200.0) (44, 43333005.1) (45, 43381785.4) (46, 43314925.5) (47, 43373693.5) (48, 43357065.4) (49, 43341847.8) (50, 43304192.3) (51, 43353422.5) (52, 43325277.3) (53, 43339375.9) (54, 43290113.7) (55, 43380262.4) (56, 43355500.8) (57, 43320207.4) (58, 43320313.4) (59, 43314854.1) (60, 43348491.9) (61, 43364193.2) (62, 43380803.3) (63, 43342847.5) (64, 43360326.9) (65, 43315938.3) (66, 43283925.0) (67, 43312418.3) (68, 43327954.0) (69, 43344369.3) (70, 43309696.6) (71, 43320083.4) (72, 43316067.0) (73, 43303258.6) (74, 43344595.3) (75, 43326390.0) (76, 43299543.9) (77, 43339194.3) (78, 43348077.9) (79, 43301171.6) (80, 43331138.3) (81, 43283485.8) (82, 43291348.6) (83, 43317357.8) (84, 43294346.3) (85, 43292978.1) (86, 43352449.8) (87, 43310141.6) (88, 43365387.8) (89, 43329165.4) (90, 43313879.1) (91, 43287974.8) (92, 43353362.9) (93, 43292705.7) (94, 43332053.3) (95, 43342681.6) (96, 43305822.3) (97, 43283503.3) (98, 43316700.3) (99, 43316082.9) (100, 43345988.5)
};
\addlegendentry{SPWL UB $|\hat{\Omega_t}| = 20$}
\label{fig:50:reg_adp_ub}
\end{axis}
\end{tikzpicture}
\caption{Numerical comparison of multistage stochastic optimization methods for $|R^x_t|= 50$.}
\label{fig:50}
\end{figure}
\begin{figure}[H]
\begin{tikzpicture}[baseline]
\pgfplotsset{every axis legend/.append style={font=\footnotesize, at={(0.6,0.0)},anchor=south, legend columns=3}}
\begin{axis}[
height=2.4in,
width=6.0in,
grid=major,
xtick = {0, 50,...,100},
ymin=41000000.0,
ymax=44500000.0,
title={Regularized SDDP vs. SPWL},
xlabel={Iteration},
ylabel={Value (\$)}]
\addplot[color=blue, dash pattern=on 2pt off 2pt, line width = 2] coordinates {
 (1, 43947539.3) (25, 43640398.1) (50, 43500826.2) (75, 43357501.4) (100, 43361648.3)
};
\addlegendentry{SDDP UB $|{\Omega_t}| = 100$}
\label{fig:100:m_sddp_ub}
\addplot[color=blue, line width = 2] coordinates {
 (1, 19318352.6) (2, 42881946.9) (3, 42882140.7) (4, 42882221.6) (5, 42882327.7) (6, 42882453.3) (7, 42882590.7) (8, 42883023.9) (9, 42883241.2) (10, 42883347.9) (11, 42883891.1) (12, 42884240.2) (13, 42884303.1) (14, 42884337.6) (15, 42884722.6) (16, 42885310.5) (17, 42885869.0) (18, 42887435.3) (19, 42888437.8) (20, 42889716.3) (21, 42891097.4) (22, 42892887.6) (23, 42894745.6) (24, 42895754.2) (25, 42896871.3) (26, 42900167.6) (27, 42900168.9) (28, 42901050.1) (29, 42905795.3) (30, 42906362.7) (31, 42906362.8) (32, 42906365.9) (33, 42906367.6) (34, 42906383.3) (35, 42906383.4) (36, 42907496.9) (37, 42909971.5) (38, 42910631.6) (39, 42911129.7) (40, 42911145.0) (41, 42911161.7) (42, 42911230.9) (43, 42911246.2) (44, 42912240.4) (45, 42912271.8) (46, 42912693.4) (47, 42912767.5) (48, 42914539.4) (49, 42914602.1) (50, 42914954.2) (51, 42915268.1) (52, 42915291.5) (53, 42915446.1) (54, 42915531.7) (55, 42915953.7) (56, 42915983.7) (57, 42916015.2) (58, 42916076.3) (59, 42916127.2) (60, 42916137.4) (61, 42916145.5) (62, 42916148.9) (63, 42916151.0) (64, 42916155.2) (65, 42916157.7) (66, 42916158.8) (67, 42916161.1) (68, 42916179.5) (69, 42916195.8) (70, 42916218.3) (71, 42916219.9) (72, 42916221.4) (73, 42916270.5) (74, 42916286.2) (75, 42916295.1) (76, 42916297.1) (77, 42916297.9) (78, 42916300.4) (79, 42916301.5) (80, 42916312.0) (81, 42916321.2) (82, 42916333.5) (83, 42916341.4) (84, 42916387.7) (85, 42916407.7) (86, 42916423.0) (87, 42916453.4) (88, 42916473.0) (89, 42916473.8) (90, 42916474.6) (91, 42916475.1) (92, 42916477.0) (93, 42919645.4) (94, 42919672.1) (95, 42919680.4) (96, 42919703.5) (97, 42919748.5) (98, 42919755.6) (99, 42919768.8) (100, 42919779.6)
};
\addlegendentry{SDDP LB $|{\Omega_t}|=100$}
\label{fig:100:m_sddp_lb}
\addplot[color=red, dash pattern=on 1pt off 1pt, line width = 2] coordinates {
 (1, 43782776.1) (25, 43418535.9) (50, 43394864.9) (75, 43357720.4) (100, 43319010.2)
};
\addlegendentry{SDDP UB $|\hat{\Omega_t}| = 20$}
\label{fig:100:m_reg_sddp_ub}
\addplot[color=red, line width = 2] coordinates {
 (1, 19255773.0) (2, 41618903.2) (3, 42440870.2) (4, 42441035.7) (5, 42606384.1) (6, 42606524.5) (7, 42688578.3) (8, 42688581.5) (9, 42688581.5) (10, 42688581.8) (11, 42688581.8) (12, 42688581.8) (13, 42688581.8) (14, 42688582.2) (15, 42688582.2) (16, 42688585.5) (17, 42688586.5) (18, 42688586.5) (19, 42688586.5) (20, 42689034.5) (21, 42689039.6) (22, 42695645.7) (23, 42695791.1) (24, 42695832.1) (25, 42695839.7) (26, 42695870.2) (27, 42697333.6) (28, 42698000.3) (29, 42698282.9) (30, 42698340.5) (31, 42698502.3) (32, 42698545.4) (33, 42698580.6) (34, 42698602.9) (35, 42698633.9) (36, 42699060.1) (37, 42699213.0) (38, 42699419.0) (39, 42699904.3) (40, 42699971.3) (41, 42700015.3) (42, 42700139.9) (43, 42700173.2) (44, 42701090.8) (45, 42701365.3) (46, 42701656.3) (47, 42701858.6) (48, 42701911.5) (49, 42702114.0) (50, 42702226.7) (51, 42702496.5) (52, 42702533.7) (53, 42702572.5) (54, 42702588.7) (55, 42703018.8) (56, 42703187.1) (57, 42703265.3) (58, 42703266.0) (59, 42703277.5) (60, 42705377.4) (61, 42751998.5) (62, 42753486.1) (63, 42753649.4) (64, 42753849.2) (65, 42753912.7) (66, 42754027.2) (67, 42754144.3) (68, 42754219.5) (69, 42754295.7) (70, 42825814.4) (71, 42835047.4) (72, 42835724.5) (73, 42836348.3) (74, 42836408.5) (75, 42836419.9) (76, 42836523.7) (77, 42845709.0) (78, 42846105.7) (79, 42846662.6) (80, 42846824.3) (81, 42848634.9) (82, 42848752.4) (83, 42848880.1) (84, 42848938.4) (85, 42848982.8) (86, 42849061.3) (87, 42849087.9) (88, 42849100.6) (89, 42849117.9) (90, 42849124.4) (91, 42849129.2) (92, 42849270.3) (93, 42850318.8) (94, 42850443.5) (95, 42850972.2) (96, 42851437.4) (97, 42851558.6) (98, 42851591.7) (99, 42851602.0) (100, 42851657.5)
};
\addlegendentry{SDDP LB $|\hat{\Omega_t}| = 20$}
\label{fig:100:m_reg_sddp_lb}
\addplot[color=green, dash pattern=on 6pt off 6pt, line width = 3] coordinates {
 (1, 44125781.6) (2, 44050313.1) (4, 44449345.5) (5, 44370696.5) (6, 44177616.7) (7, 44023161.0) (8, 43953791.1) (9, 43889516.5) (10, 43891714.0) (11, 43840984.5) (12, 43776789.1) (13, 43777769.7) (14, 43767262.1) (15, 43759718.5) (16, 43751687.3) (17, 43729603.1) (18, 43695695.3) (19, 43717580.5) (20, 43702468.7) (21, 43762274.5) (22, 43648488.6) (23, 43715936.4) (24, 43623274.2) (25, 43720096.2) (26, 43647260.9) (27, 43622814.2) (28, 43645827.7) (29, 43752959.9) (30, 43616112.4) (31, 43653658.4) (32, 43544678.9) (33, 43572921.9) (34, 43549303.4) (35, 43520032.2) (36, 43558807.6) (37, 43528488.1) (38, 43567644.5) (39, 43559712.1) (40, 43500489.7) (41, 43492544.5) (42, 43520494.7) (43, 43531649.3) (44, 43496207.8) (45, 43475172.9) (46, 43477639.1) (47, 43485709.5) (48, 43482495.5) (49, 43486722.5) (50, 43450632.8) (51, 43505027.2) (52, 43432360.9) (53, 43448414.2) (54, 43442280.6) (55, 43467693.7) (56, 43416753.2) (57, 43410865.2) (58, 43479311.7) (59, 43408592.5) (60, 43401124.5) (61, 43409761.1) (62, 43430217.1) (63, 43437625.3) (64, 43383081.1) (65, 43475480.5) (66, 43463899.2) (67, 43381425.4) (68, 43444899.5) (69, 43396101.8) (70, 43421118.8) (71, 43382889.6) (72, 43461355.0) (73, 43396698.1) (74, 43402184.2) (75, 43328162.8) (76, 43331803.9) (77, 43393561.6) (78, 43426405.5) (79, 43371495.5) (80, 43398384.4) (81, 43405916.3) (82, 43387396.2) (83, 43368163.7) (84, 43355634.1) (85, 43378748.9) (86, 43424396.9) (87, 43384269.3) (88, 43344763.2) (89, 43354512.5) (90, 43385082.8) (91, 43362136.2) (92, 43349923.0) (93, 43334874.7) (94, 43340453.1) (95, 43368557.3) (96, 43310114.0) (97, 43349168.0) (98, 43346886.6) (99, 43358735.0) (100, 43335658.6)
};
\addlegendentry{SPWL UB $|{\Omega_t}| = 100$}
\label{fig:100:reg_adp_ub}
\addplot[color=orange, dash pattern=on 4pt off 4pt, line width = 2] coordinates {
 (1, 43786418.6) (2, 43853061.5) (3, 44221103.5) (4, 44140392.0) (5, 43832281.4) (6, 43762075.0) (7, 43705194.6) (8, 43653332.2) (9, 43607862.9) (10, 43579924.9) (11, 43632538.1) (12, 43605419.0) (13, 43586635.3) (14, 43624907.2) (15, 43560825.5) (16, 43564138.7) (17, 43592077.7) (18, 43528912.8) (19, 43524535.9) (20, 43501709.8) (21, 43524766.1) (22, 43496536.9) (23, 43486830.2) (24, 43480722.3) (25, 43518476.2) (26, 43483815.2) (27, 43483105.1) (28, 43453007.8) (29, 43470679.7) (30, 43442740.4) (31, 43440746.3) (32, 43421336.6) (33, 43450895.5) (34, 43406627.1) (35, 43412696.0) (36, 43437560.5) (37, 43433397.0) (38, 43390939.0) (39, 43412914.8) (40, 43375237.3) (41, 43385620.7) (42, 43377547.2) (43, 43383534.2) (44, 43392337.4) (45, 43422236.8) (46, 43352414.2) (47, 43394813.2) (48, 43327271.5) (49, 43343392.1) (50, 43325130.1) (51, 43411369.3) (52, 43387274.3) (53, 43365223.4) (54, 43343963.8) (55, 43344920.7) (56, 43333088.7) (57, 43345445.1) (58, 43309112.8) (59, 43313049.8) (60, 43316366.2) (61, 43369196.5) (62, 43327017.2) (63, 43326143.7) (64, 43378633.1) (65, 43302597.8) (66, 43281692.2) (67, 43283556.8) (68, 43362405.7) (69, 43336898.2) (70, 43288732.6) (71, 43284480.1) (72, 43297790.3) (73, 43296853.7) (74, 43308819.7) (75, 43330313.0) (76, 43296155.8) (77, 43306424.8) (78, 43274809.0) (79, 43264879.2) (80, 43326307.7) (81, 43298596.6) (82, 43257118.0) (83, 43234369.5) (84, 43300709.3) (85, 43303811.1) (86, 43291945.2) (87, 43265457.0) (88, 43288877.5) (89, 43271321.8) (90, 43269891.7) (91, 43308765.5) (92, 43285447.9) (93, 43268291.5) (94, 43240426.7) (95, 43251667.7) (96, 43281447.7) (97, 43287531.3) (98, 43265135.4) (99, 43256275.7) (100, 43257721.9)
};
\addlegendentry{SPWL UB $|\hat{\Omega_t}| = 20$}
\label{fig:100:reg_adp_ub}
\end{axis}
\end{tikzpicture}
\caption{Numerical comparison of multistage stochastic optimization methods for $|R^x_t|= 100$.}
\label{fig:100}
\end{figure}
% \normalsize
\begin{table}
\normalsize
\begin{center}
\begin{tabular}{|l|c|c|c|c|}
\hline
\multicolumn{2}{|c|}{\backslashbox[13em]{$(|R^x_t|, |\hat{\Omega}_t|)$ \qquad Method}{\# Iterations}} & 1 & 50 & 100 \\
\hline
(25, 50) 
& \begin{tabular}{@{}c@{}}SDDP (Regularization) \\ADP - SPWL\end{tabular} 
&  \begin{tabular}{@{}c@{}}712.0\\ 284.0\end{tabular}
&  \begin{tabular}{@{}c@{}}804.4\\ 310.2\end{tabular}
&  \begin{tabular}{@{}c@{}}902.1\\ 325.0\end{tabular}
\\ \hline
(25, 20) 
& \begin{tabular}{@{}c@{}}SDDP (Regularization) \\ADP - SPLW\end{tabular} 
&  \begin{tabular}{@{}c@{}}176.0\\ 116.0\end{tabular}
&  \begin{tabular}{@{}c@{}}180.7\\ 136.4\end{tabular}
&  \begin{tabular}{@{}c@{}}185.4\\ 147.1\end{tabular}
\\ \hline
(50, 50) 
& \begin{tabular}{@{}c@{}}SDDP (Regularization) \\ADP - SPWL\end{tabular} 
&  \begin{tabular}{@{}c@{}}763.0\\ 719.0\end{tabular}
&  \begin{tabular}{@{}c@{}}869.1\\ 762.1\end{tabular}
&  \begin{tabular}{@{}c@{}}986.4\\ 795.8\end{tabular}
\\ \hline
(50, 50) 
& \begin{tabular}{@{}c@{}}SDDP (Regularization) \\ADP - SPWL\end{tabular} 
&  \begin{tabular}{@{}c@{}}451.0\\ 975.0\end{tabular}
&  \begin{tabular}{@{}c@{}}487.4\\ 1081.9\end{tabular}
&  \begin{tabular}{@{}c@{}}517.6\\ 964.5\end{tabular}
\\ \hline
(100, 100) 
& \begin{tabular}{@{}c@{}}SDDP (Regularization) \\ADP - SPWL\end{tabular} 
&  \begin{tabular}{@{}c@{}}2229.0\\ 3167.0\end{tabular}
&  \begin{tabular}{@{}c@{}}2657.1\\ 3016.0\end{tabular}
&  \begin{tabular}{@{}c@{}}3030.8\\ 2992.7\end{tabular}
\\ \hline
(100, 20) 
& \begin{tabular}{@{}c@{}}SDDP (Regularization) \\ADP - SPWL\end{tabular} 
&  \begin{tabular}{@{}c@{}}478.0\\ 3167.0\end{tabular}
&  \begin{tabular}{@{}c@{}}616.9\\ 3016.0\end{tabular}
&  \begin{tabular}{@{}c@{}}683.9\\ 2992.7\end{tabular}
\\ \hline
\end{tabular}
\caption{Computational time per iteration (in seconds) for stochastic optimization methods.}
\label{t:swi}
\end{center}
\end{table}
\normalsize

Table \ref{t:swi} shows the CPU times (in seconds) per iteration for problems with up to 100 storage devices, and different sample sizes $\hat{\Omega}_t$.
We can see that computational times of all methods depend on the choice of problem parameters, and there is not a a single method that always outperforms the competition.
Hence, the choice of the solution method should be made on a case-by-case basis taking into account factors such as prefered solution structure, limitations on computational time and resources, ease of implementation, the characteristics of the stochastic process, the quality of the available linear and quadratic programming solvers, and potential future needs.

In addition, we are also interested in the magnitude of the errors arising from the use of a sampled model. In order to gain estimates, we consider runs with a sample sizes $|\Omega_t| = 100$, and then do smaller samples $\hat{\Omega_t}$ drawn from $\Omega_t$, and compare the resulting objective values.
The use of a small sample $\hat{\Omega}_t$ could lead to a major computational speed-up since SDDP-type methods need to solve an optimization problem for each random realization at each time step in the backward pass of every iteration.
However, we can also see that the use of a small sample $\hat{\Omega}_t$ can also introduce significant approximation errors. It is widely believed (and theoretically expected) that a smaller sample size should result in a lower optimal objective value due to overfitting. However, we can see that in practice that is not necessarily the case since we rarely solve the problem to optimality. Instead a smaller sample can result in both underestimated, as well as overestimated bounds and there is no obvious criterion to
help us distinguish between the two cases a priori.
In addition, it is known that the square root law applies to the optimal values of the sampled optimization problems (i.e. in order to gain one decimal place of precision, we need to increase the sample size by a factor of 100) but large-scale real world problems are rarely solved to optimality. 
% In addition, for a given number of iterations.
\input{temp4.tex}
\input{temp5.tex}
\input{temp6.tex}
Moreover, practitioners usually use a small fixed number of iterations (hundreds or thousands) and hence there is no direct way to determine how a smaller sample size would affect the resulting policy and optimality gap.  
Finally, we would like to emphasize that typically sampling is not an issue for approximate dynamic programming with separable value functions since it works directly with the full probability model, and hence it avoids sampling errors altogether.
Unfortunately, currently there does not exist any approach that would allow us to extend this property to the SDDP framework.
% Still, we include runs of ADP-SPWL with different sets $\hat{\Omega}_t$ for the sake of comparison. 

\section{Conclusion}
\label{sec:gridpaperconclusions}

Multistage stochastic optimization problems with long time horizons appear in various fields and the computational solution of such models is a topic of growing importance.
In our work we have compared the performance of SDDP (with regularization) and ADP-SPWL for the solution of multistage stochastic problems in energy storage.
% , and discussed their characteristics.

On one hand, approximations with separable value functions feature many desirable properties.
For instance, the random process can depend on the decisions made at each time step, and use of a wide variety of probability models is possible (even computer simulations).
However, such broad applicability comes at a price. In general, practitioners do not have any guarantees for the quality of policies derived with separable value functions since lower bounds are not readily available.

When optimality guarantees are needed, researchers often resort to the SDDP framework which requires a memoryless stochastic process that is exogenous to the decisions made by the model.
Moreover, SDDP employs a finite sample to construct a sampled average approximation to the original problem. When the given formulation is convex, lower bounds are readily available and an optimality bound for the policy of the sampled problem can be estimated.
However, our numerical work indicates that the relationship between the bounds derived from the sampled model, and the bounds derived from the full model cannot be determined a priori.  

In general, the practical application of SDDP has been limited to problems
with low--dimensional value functions such as hydro--power problems with a small number of (groups of) reservoirs. In our experiments we show that such a limitation can be overcome with the use of a regularization technique proposed by \cite{asamov2015regularized}, and approximations with Benders cuts can compete with ADP methods with separable value functions on larger problems than previously known. Moreover, we have studied the resulting policies and compared them to the optimal solutions in deterministic instances. Our results suggest that neither algorithm is universally best, and we recommend that practitioners choose the solution method depending on the characteristics of the problem at hand. 

In the future we plan to extend the current work to numerical testing of non-linear problems, including risk--averse models formulated as time--consistent compositions of coherent measures of risk (\cite{asamov2014risk}), as well as problems with multiple objectives (\cite{young2010multiobjective}).  
%Additionally, we could also explore possible applications for the solution of multistage mixed-integer stochastic models.
%

\bibliographystyle{mynatbib} % (uses file "plain.bst")
\bibliography{myrefs,library}   % expects file "myrefs.bib"

% Acknowledgments here
% \ACKNOWLEDGMENT{%
% % Enter the text of acknowledgments here
% }% Leave this (end of acknowledgment)

% Appendix here
% Options are (1) APPENDIX (with or without general title) or 
%             (2) APPENDICES (if it has more than one unrelated sections)
% Outcomment the appropriate case if necessary
%
%
%   or 
%
% \begin{APPENDICES}
% \section{<Title of Section A>}
% \section{<Title of Section B>}
% etc
% \end{APPENDICES}

% References here (outcomment the appropriate case) 

% CASE 1: BiBTeX used to constantly update the references 
%   (while the paper is being written).
%\bibliographystyle{informs2014} % outcomment this and next line in Case 1
%\bibliography{<your bib file(s)>} % if more than one, comma separated

% CASE 2: BiBTeX used to generate mypaper.bbl (to be further fine tuned)
%\input{mypaper.bbl} % outcomment this line in Case 2

\end{document}